\documentclass[a4paper,11pt]{amsart}
\usepackage{amssymb,amsfonts,amsxtra,amscd}
\usepackage[all]{xy}
\usepackage{fullpage}
\theoremstyle{plain}
\newtheorem{theorem}{Theorem}[section]
\newtheorem{lemma}[theorem]{Lemma}

\newtheorem{prop}[theorem]{Proposition}
\theoremstyle{definition}
\newtheorem{defi}[theorem]{Definition}

\theoremstyle{remark}
\newtheorem{rem}[theorem]{Remark}
\numberwithin{equation}{section}
\newcommand{\ci}{\ensuremath{C_\infty}}
\newcommand{\ai}{\ensuremath{A_\infty}}
\newcommand{\li}{\ensuremath{L_\infty}}
\newcommand{\clalg}[1]{\ensuremath{\widehat{L}\Sigma #1^*}}
\newcommand{\deof}[1]{\ensuremath{\Omega^1(#1)}}
\newcommand{\drof}[1]{\ensuremath{DR^1(#1)}}
\newcommand{\drnof}[1]{\ensuremath{\overline{DR}^1(#1)}}
\newcommand{\drzf}[1]{\ensuremath{DR^0(#1)}}
\newcommand{\drnzf}[1]{\ensuremath{\overline{DR}^0(#1)}}
\newcommand{\drtf}[1]{\ensuremath{DR^2(#1)}}
\newcommand{\de}[1]{\ensuremath{\Omega^\bullet(#1)}}
\newcommand{\dr}[1]{\ensuremath{DR^\bullet(#1)}}
\newcommand{\mstr}[4][]{\ensuremath{\mathcal{#1#2}_{#3}(#4)}}
\newcommand{\mext}[3][]{\ensuremath{\mathcal{#1E}_{#2}(#3)}}
\newcommand{\mmor}[4][]{\ensuremath{\mathcal{#1M}_{#2}(#3 ; #4)}}
\newcommand{\mmext}[5][]{\ensuremath{\mathcal{#1EM}_{#2}(#3 : #4 , #5)}}
\newcommand{\caq}[3][\bullet]{\ensuremath{C^{#1}_\mathrm{Harr}(#2,#3)}}
\newcommand{\haq}[3][\bullet]{\ensuremath{H^{#1}_\mathrm{Harr}(#2,#3)}}
\newcommand{\cnaq}[3][\bullet]{\ensuremath{\overline{C}^{#1}_\mathrm{Harr}(#2,#3)}}
\newcommand{\Scnaq}[3][\bullet]{\ensuremath{\overline{SC}^{#1}_\mathrm{Harr}(#2,#3)}}

\newcommand{\ccaq}[2][\bullet]{\ensuremath{CC^{#1}_\mathrm{Harr}(#2)}}
\newcommand{\hcaq}[2][\bullet]{\ensuremath{HC^{#1}_\mathrm{Harr}(#2)}}
\newcommand{\ccnaq}[2][\bullet]{\ensuremath{\overline{CC}^{#1}_\mathrm{Harr}(#2)}}
\newcommand{\hcnaq}[2][\bullet]{\ensuremath{\overline{HC}^{#1}_\mathrm{Harr}(#2)}}
\newcommand{\hhoch}[3][\bullet]{\ensuremath{H^{#1}_\mathrm{Hoch}(#2,#3)}}
\newcommand{\hchoch}[2][\bullet]{\ensuremath{HC^{#1}_\mathrm{Hoch}(#2)}}
\newcommand{\cotimes}{\ensuremath{\hat{\otimes}}}
\newcommand{\gf}{\ensuremath{\mathbb{K}}}
\newcommand{\innprod}[1]{\ensuremath{\langle -,- \rangle:#1 \otimes #1 \to \gf}}
\newcommand{\noproof}{\begin{flushright} \ensuremath{\square} \end{flushright}}
\DeclareMathOperator{\Hom}{Hom}
\DeclareMathOperator{\Der}{Der}
\DeclareMathOperator{\Aut}{Aut}
\DeclareMathOperator{\Obs}{Obs}
\DeclareMathOperator{\obs}{obs}

\DeclareMathOperator{\ad}{ad}
\DeclareMathOperator{\id}{id}

\begin{document}
\begin{abstract}
In this paper we show that a strongly homotopy commutative (or $C_\infty$-) algebra with an invariant inner product on its cohomology can be uniquely extended to a symplectic $C_\infty$-algebra (an $\infty$-generalisation of a commutative Frobenius algebra introduced by Kontsevich). This result relies on the algebraic Hodge decomposition of the cyclic Hochschild cohomology of a $\ci$-algebra and does not generalize to algebras over other operads.
\end{abstract}
\title{Symplectic $\ci$-algebras}
\author{Alastair Hamilton \and Andrey Lazarev}
\address{Mathematics Department, University of Leicester, Leicester, England. LE1 7RH.} \email{hamilton@mpim-bonn.mpg.de \and al179@le.ac.uk}
\maketitle
\tableofcontents

\section{Introduction}

The notion of a strongly homotopy associative algebra or an $A_\infty$-algebra was introduced in \cite{staha2} for the purposes of studying $H$-spaces. It was originally defined via a system of higher multiplication maps satisfying a series of complicated relations. It is an algebra over the operad of planar trees; other operads give rise to various other flavours of strongly homotopy algebras, cf. \cite{markl}. In this paper we are concerned chiefly with the homotopy invariant generalisation of commutative algebras, the so-called $\ci$-algebras. More precisely, we shall study \emph{symplectic} $\ci$-algebras; this is a strong homotopy analogue of a commutative Frobenius algebra. Symplectic infinity-algebras (also called cyclic infinity algebras or infinity algebras with an invariant inner product) were introduced in \cite{kontfd} and \cite{kontsg} and were shown to have a close relation with graph homology and therefore to the intersection theory on the moduli spaces of complex curves and invariants of differentiable manifolds.  A short, informal introduction to graph homology is contained in \cite{voronov}; a more substantial account is in \cite{vogtmann}. The connections between graph complexes and infinity algebras are studied in \cite{hamgraph} and \cite{ccs}.

Our main theorem here states a unital $\ci$-algebra with the structure of a Frobenius algebra on its cohomology is weakly equivalent to a symplectic $\ci$-algebra which is uniquely determined up to homotopy. This is a rather surprising phenomenon which does not hold in the context of other infinity-algebras. One application of this result is that the cohomology algebra of a simply-connected manifold supports the structure of a symplectic $C_\infty$-algebra which leads to a manifestly homotopy invariant construction of  string topology-type operations, cf. \cite{string1, string2}. Another potential application is to the study of the derived category of coherent sheaves on a complex manifold; note that according to \cite{block} the latter could be realised as a certain localisation of the category of differential-graded modules over the Dolbeault algebra, which supports an invariant inner product in the Calabi-Yau case.

The basic fact which leads to this result is that the cohomology theory for $\ci$-algebras, known as Harrison cohomology, essentially coincides with the cohomology theory for symplectic $\ci$-algebras, called cyclic Harrison cohomology. This result in turn follows from the algebraic Hodge decomposition for the cyclic Hochschild cohomology of a $\ci$-algebra and is proved in the author's paper \cite{HL}. We shall use the results of this paper extensively.

To prove the main theorem we develop an obstruction theory for lifting the  symplectic $C_n$-structures and $C_n$-morphisms which is of independent interest. The proof is carried out by combining this obstruction theory with the Hodge decomposition for cyclic cohomology.

The paper is organised as follows. In Section \ref{sec_noncom}  we introduce the basics of noncommutative geometry. Section \ref{sec_calg} discusses $\ci$-algebras, with and without an invariant inner product, and the associated cohomology theories. These sections are merely recollections of the main results of \cite{HL} and contain almost no proofs.

Section \ref{sec_obstrc} deals with the obstruction theory for $C_n$-structures and morphisms, in both the symplectic and the nonsymplectic contexts. In Section \ref{correspondence} we prove the theorem mentioned above on the relationship between symplectic and nonsymplectic $\ci$-algebras. The cases of unital and non-unital symplectic $\ci$-algebras differ slightly and we treat each of them in their own detail.

\subsection{Notation and conventions} \label{sec_notcon}

We adopt the notation and conventions of \cite{HL}. Throughout the paper our ground ring $\gf$ will be an evenly graded commutative ring containing the
field $\mathbb{Q}$.  \gf-algebras and \gf-modules will simply be called algebras and modules. We will assume that all of our \gf-modules are obtained from $k$-vector spaces by extension of scalars where $k=\mathbb{Q}$ or any other subfield of $\gf$. All of our tensors will be taken over the ground ring {\gf} unless stated otherwise.

Given a graded module V we define the tensor algebra $TV$ by
\[ TV:=\gf \oplus V \oplus V^{\otimes 2} \oplus \ldots \oplus V^{\otimes n} \oplus \ldots. \]
We define the free Lie algebra $LV$ as the Lie subalgebra of the commutator algebra $TV$ which consists of linear combinations of Lie monomials.

We will use $~\widehat{ \ }~$ to denote completion. Given a profinite graded module $V$ we can define the completed (or pro-free) versions  $TV$ and $LV$. For instance the pro-free associative algebra would be
\[ \widehat{T}V :=\prod_{i=0}^{\infty} V^{\cotimes i} \]
where $\cotimes$ denotes the completed tensor product. The pro-free Lie algebra  $\widehat{L}V$ is the Lie subalgebra of $\widehat{T}V$ consisting of all convergent (possibly uncountably) infinite linear combinations of Lie monomials.

Given a profinite graded module $V$, the Lie algebra consisting of all \emph{continuous} derivations
\[ \xi : \widehat{L}V \to \widehat{L}V \]
will be denoted by $\Der(\widehat{L}V)$. In order to emphasise our use of geometrical ideas in this paper we will use the term `vector field' synonymously with `continuous derivation'. The group consisting of all invertible \emph{continuous} Lie algebra homomorphisms
\[ \phi: \widehat{L}V \to \widehat{L}V \]
will be denoted by $\Aut(\widehat{L}V)$. Again, in order to emphasise the geometrical approach we will call an `invertible continuous homomorphism' of formal graded commutative, associative or Lie algebras a `diffeomorphism'.

Given a profinite graded module $V$ we can place a grading on $\widehat{T}V$ which is different from the grading which is naturally inherited from $V$. We say an element $x \in \widehat{T}V$ has homogeneous \emph{order} $n$ if $x \in V^{\cotimes n}$. This grading naturally descends to the submodule $\widehat{L}V$.

We say a continuous endomorphism (linear map) $f:\widehat{L}V \to \widehat{L}V$ has homogeneous \emph{order} $n$ if it takes any element of order $i$ to an element of order $i+n-1$. Any vector field $\xi \in \Der(\widehat{L}V)$ could be written in the form
\begin{equation} \label{eqn_vecfrm}
\xi = \xi_1 + \xi_2 + \xi_3 + \ldots + \xi_n + \ldots,
\end{equation}
where $\xi_i$ is a vector field of order $i$.

Likewise any diffeomorphism $\phi \in \Aut(\widehat{L}V)$ could be written in the form
\[ \phi = \phi_1 + \phi_2 + \phi_3 + \ldots + \phi_n + \ldots, \]
where $\phi_i$ is an endomorphism of order $i$. We call $\phi$ a \emph{pointed diffeomorphism} if $\phi_1 = \id$.

Our convention will be to always work with cohomologically graded objects and we consequently define the suspension of a graded module $V$ as $\Sigma V$ where $\Sigma V^i:=V^{i+1}$. We define the desuspension of $V$ as $\Sigma^{-1}V$ where $\Sigma^{-1}V^{i}:= V^{i-1}$. The term `differential graded algebra' will be abbreviated as `DGA'.

Given a graded module $V$, we denote the graded \gf-linear dual
\[ \Hom_\gf(V,\gf) \]
by $V^*$. For the sake of clarity, when we write $\Sigma V^*$ we mean the graded module $\Hom_\gf(\Sigma V,\gf)$.

\section{Noncommutative geometry}\label{sec_noncom}

\subsection{Formal noncommutative geometry}

This section reviews the basics of formal noncommutative Lie calculus including the de Rham complex, the Lie derivative along a vector field and related constructions.

\begin{defi}
Let $W$ be a free graded module.
Consider the module $\widehat{T}(W^*) \cotimes \widehat{L}(W^*)$ and write $x\cotimes y=x\cotimes dy$. Then the module of Lie 1-forms $\deof{W}$ is defined as the quotient of $\widehat{T}(W^*) \cotimes \widehat{L}(W^*)$ by the relations
\[ x \cotimes d[y,z]=xy \cotimes dz -(-1)^{|z||y|}xz \cotimes dy. \]
This is a left $\widehat{L}(W^*)$-module via the action
\[ x \cdot (y \cotimes dz):= xy \cotimes dz. \]
\end{defi}

\begin{defi}
Let $W$ be a free graded module. The module of Lie forms $\de{W}$ is defined as
\[ \de{W}:=\widehat{L}(W^*)\ltimes \widehat{L}(\Sigma^{-1}\deof{W}) \]
where the action of $\widehat{L}(W^*)$ on $\widehat{L}(\Sigma^{-1}\deof{W})$ is the restriction of the standard action of $\widehat{L}(W^*)$ on $\widehat{T}(\Sigma^{-1}\deof{W})$ to the Lie subalgebra of Lie monomials. The map $d:\widehat{L}(W^*) \to \deof{W}$ lifts uniquely to give $\de{W}$ the structure of a formal DGLA.
\end{defi}

\begin{rem} \label{rem_form_func}
This construction is functorial; i.e. given free graded modules $V$ and $W$ and a homomorphism of Lie algebras $\phi:\widehat{L}(V^*) \to \widehat{L}(W^*)$ (in particular, any linear map from $W$ to $V$ gives rise to such a map), there is a unique map $\phi^*:\de{V} \to \de{W}$ extending $\phi$ to a homomorphism of differential graded Lie algebras.
\end{rem}

\begin{defi} \label{def_operators}
Let $W$ be a free graded module and let $\xi:\widehat{L}(W^*) \to \widehat{L}(W^*)$ be a vector field:
\begin{enumerate}
\item[(i)]
We can define a vector field $L_\xi: \de{W} \to \de{W}$, called the Lie derivative, by the formulae;
\begin{displaymath}
\begin{array}{lc}
L_\xi(x):=\xi(x), & x \in \widehat{L}(W^*); \\
L_\xi(dx):=(-1)^{|\xi|}d(\xi(x)), & x \in \widehat{L}(W^*). \\
\end{array}
\end{displaymath}
\item[(ii)]
We can define a vector field $i_\xi:\de{W} \to \de{W}$ of bidegree $(-1,|\xi|)$, called the contraction operator, by the formulae;
\begin{displaymath}
\begin{array}{lc}
i_\xi(x):=0, & x \in \widehat{L}(W^*); \\
i_\xi(dx):=\xi(x), & x \in \widehat{L}(W^*). \\
\end{array}
\end{displaymath}
\end{enumerate}
\end{defi}

These operators satisfy certain important identities which are summarised by the following lemma.

\begin{lemma} \label{lem_schf}
Let $V$ and $W$ be free graded modules. Let $\xi,\gamma:\widehat{L}(V^*) \to \widehat{L}(V^*)$ be vector fields and let $\phi:\widehat{L}(V^*) \to \widehat{L}(W^*)$ be a diffeomorphism, then we have the following identities:
\begin{displaymath}
\begin{array}{rl}
\textnormal{(i)} & L_\xi=[i_\xi,d]. \\
\textnormal{(ii)} & [L_\xi,i_\gamma]=i_{[\xi,\gamma]}. \\
\textnormal{(iii)} & L_{[\xi,\gamma]}=[L_\xi,L_\gamma]. \\
\textnormal{(iv)} & [i_\xi,i_\gamma]=0. \\
\textnormal{(v)} & [L_\xi,d]=0. \\
\textnormal{(vi)} & L_{\phi\xi\phi^{-1}}=\phi^* L_\xi \phi^{*-1}. \\
\textnormal{(vii)} & i_{\phi\xi\phi^{-1}}=\phi^* i_\xi \phi^{*-1}. \\
\end{array}
\end{displaymath}
\end{lemma}
\noproof

\begin{defi}
Let $W$ be a free graded module.
The de Rham complex $\dr{W}$ is defined as the quotient of $\de{W} \cotimes \de{W}$ by the relations
\begin{displaymath}
\begin{array}{ll}
x \cotimes y = (-1)^{|x||y|}y \cotimes x; & x,y \in \de{W}; \\
x \cotimes [y,z] = [x,y] \cotimes z; & x,y,z \in \de{W}. \\
\end{array}
\end{displaymath}
The differential $d:\dr{W} \to \dr{W}$ is induced by the differential $d:\de{W} \to \de{W}$ by specifying
\[ d(x \cotimes y):=dx \cotimes y + (-1)^{|x|}x \cotimes dy; \quad x,y \in \de{W}. \]
\end{defi}

\begin{rem}
It follows from Remark \ref{rem_form_func} that the construction of the de Rham complex is also functorial; that is to say that if $\phi:\widehat{L}(V^*) \to \widehat{L}(W^*)$ is a homomorphism of algebras then it induces a morphism of complexes $\phi^*:\dr{V} \to \dr{W}$ defined by the formula
\[ \phi^*(x \cotimes y):= \phi^*(x) \cotimes \phi^*(y); \quad x,y \in \de{V}. \]

Furthermore the Lie and contraction operators $L_\xi,i_\xi:\de{V}\to\de{V}$ defined by Definition \ref{def_operators} factor naturally through Lie and contraction operators $L_\xi,i_\xi:\dr{V}\to\dr{V}$, which by an abuse of notation, we denote by the same letters. The Lie operator is given by the formula
\[ L_\xi(x \cotimes y):= L_\xi(x) \cotimes y + (-1)^{|x|}x \cotimes L_\xi(y); \quad x,y \in \de{V}. \]
The contraction operator is defined in the same way.
\end{rem}

\subsection{Formal noncommutative symplectic geometry}
We will begin by recalling the basic terminology of symplectic geometry.

\begin{defi} \label{def_sympform}
Let $W$ be a free graded module. A homogeneous 2-form $\omega \in \drtf{W}$ is nondegenerate if and only if the map $\Phi:\Der(\widehat{L}(W^*)) \to \drof{W}$ of degree $|\omega|$ defined by the formula
\begin{equation} \label{eqn_ofmvec}
\Phi(\xi):= i_\xi(\omega), \quad \xi \in \Der(\widehat{L}(W^*))
\end{equation}
is bijective.

Furthermore, we say that $\omega$ is a symplectic form if it is also closed (i.e. $dw=0$).
\end{defi}

We say a 2-form $\omega\in\drtf{W}$ is a  \emph{constant} 2-form if it can be written in the form $\omega=\sum_i dy_i dz_i$ for some functions $y_i,z_i\in W^*$.

\begin{prop}
Let $W$ be a free graded module. There exists a map
\[ \kappa:\{ \omega \in \drtf{W} : \omega \text{ is constant} \} \to (\Lambda^2 W)^* \]
defined by the formula
\[ \kappa(dxdy):=(-1)^x[x\otimes y-(-1)^{xy}y\otimes x] \]
which provides an isomorphism between constant 2-forms and \emph{skew-symmetric} bilinear forms on $W$. Furthermore, a constant 2-form $\omega\in\drtf{W}$ is nondegenerate if and only if the corresponding bilinear form $\langle-,-\rangle:=\kappa(\omega)$ is nondegenerate.
\end{prop}
\noproof

\begin{defi} \label{def_symvec}
Let $V$ and $W$ be free graded modules and let $\omega \in \drtf{V}$ and $\omega' \in \drtf{W}$ be homogeneous symplectic forms:
\begin{enumerate}
\item[(i)]
We say a vector field $\xi:\widehat{L}(V^*) \to \widehat{L}(V^*)$ is a \emph{symplectic vector field} if $L_\xi(\omega)=0$.
\item[(ii)]
We say a diffeomorphism $\phi:\widehat{L}(V^*) \to \widehat{L}(W^*)$ is a \emph{symplectomorphism} if $\phi^*(\omega)=\omega'$.
\end{enumerate}
\end{defi}

We have the following simple proposition about symplectic vector fields which can be found in \cite{ginzsg}:

\begin{prop} \label{prop_symofm}
Let $W$ be a free graded module and let $\omega \in \drtf{W}$ be a symplectic form. Then the map $\Phi$ defined by equation \eqref{eqn_ofmvec} induces a one-to-one correspondence between \emph{symplectic} vector fields and \emph{closed} 1-forms:
\[ \Phi: \{ \xi \in \Der(\widehat{L}(W^*)) : L_\xi(\omega) = 0 \} \to \{ \alpha \in \drof{W} : d\alpha = 0 \}. \]
\end{prop}
\noproof

\section{$\ci$-algebras and related cohomology theories}\label{sec_calg}

In this section we consider the notion of a $\ci$-algebra and that of a symplectic $C_\infty$-algebra. Our presentation will be based on the noncommutative Lie geometry reviewed in the previous section.

\subsection{Main definitions}

\begin{defi} \label{def_infstr}
Let $V$ be a free graded module. A $\ci$-structure on $V$ is a vector field
\[ m:\clalg{V} \to \clalg{V} \]
of degree one, such that $m^2=0$.
\end{defi}

\begin{defi} \label{def_infmor}
Let $V$ and $U$ be free graded modules. Let $m$ and $m'$ be $\ci$-structures on $V$ and $U$ respectively. A $\ci$-morphism from $V$ to $U$ is a continuous algebra homomorphism
\[ \phi:\clalg{U} \to \clalg{V} \]
of degree zero such that $\phi \circ m'=m \circ \phi$.
\end{defi}

\begin{defi}
Let $V$ be a free graded module:
\begin{enumerate}
\item
We say that a $\ci$-structure $m:\clalg{V} \to \clalg{V}$ is unital if there is a distinguished element $1 \in V$ (the unit) of degree zero which can be extended to a basis $1,\{x_i\}_{i\in I}$ of $V$ such that $m$ has the form;
\begin{equation} \label{eqn_unitform}
m=A(\boldsymbol{t})\partial_\tau+\sum_{i\in I}B_i(\boldsymbol{t})\partial_{t_i} + \ad\tau - \tau^2\partial_\tau,
\end{equation}
where $\tau,\boldsymbol{t}:=\tau,\{t_i\}_{i\in I}$ is the topological basis of $\Sigma V^*$ which is dual to the basis $\Sigma 1,\{\Sigma x_i\}_{i\in I}$ of $\Sigma V$ and $A(\boldsymbol{t}),\{B_i(\boldsymbol{t})\}_{i\in I}$ are Lie power series in $\boldsymbol{t}$. In this case we say that the $\ci$-algebra $V$ is unital or that it has a unit.
\item
Suppose that $V$ and $U$ are two unital $\ci$-algebras. We say that a $\ci$-morphism $\phi:\clalg{U} \to \clalg{V}$ is unital if $\phi$ has the form;
\begin{displaymath}
\begin{array}{ccc}
\phi(\tau') & = & \tau + A(\boldsymbol{t}), \\
\phi(t'_i) & = & B_i(\boldsymbol{t}); \\
\end{array}
\end{displaymath}
where $\tau,\boldsymbol{t}$ and $\tau',\boldsymbol{t}'$ are the topological bases of $\Sigma V^*$ and $\Sigma U^*$ which are dual to the bases $\Sigma 1_V,\{\Sigma x_i\}_{i\in I}$ and $\Sigma 1_U,\{\Sigma x'_j\}_{j\in J}$ of $\Sigma V$ and $\Sigma U$ respectively.
\end{enumerate}
\end{defi}

\subsection{Cohomology theories}

There are cohomology theories associated to $\ci$-algebras which are the home for various obstruction groups. These theories are called the \emph{Harrison cohomology} theories. In the case of a \emph{unital} $\ci$-algebra, the relevant cohomology theory is called \emph{normalised} Harrison cohomology.

\subsubsection{Unnormalised theories}

\begin{defi}
Let $V$ be a free graded module. Let $m:\clalg{V} \to \clalg{V}$ be a $\ci$-structure. The Harrison complex of the $\ci$-algebra $V$ with coefficients in $V$ is defined on the module consisting of all vector fields on $\clalg{V}$:
\[ \caq{V}{V}:=\Sigma^{-1}\Der(\clalg{V}). \]
The differential $d:\caq{V}{V} \to \caq{V}{V}$ is given by
\[ d(\xi):=[m,\xi], \quad \xi \in \Der(\clalg{V}). \]
The Harrison cohomology of $V$ with coefficients in $V$ is defined as the cohomology of the complex $\caq{V}{V}$ and denoted by $\haq{V}{V}$.
\end{defi}

\begin{rem}
The purpose of the desuspensions in the above definition is to make the grading consistent with the classical grading on these cohomology theories.
\end{rem}

\begin{defi}
Let $V$ be a free graded module and $m:\clalg{V} \to \clalg{V}$ be a $\ci$-structure. The Harrison complex of the $\ci$-algebra $V$ with coefficients in $V^*$ is defined on the module consisting of all 1-forms:
\[ \caq{V}{V^*}:=\Sigma\drof{\Sigma V}. \]
The differential on this complex is the (suspension of the) Lie operator of the vector field $m$;
\[ L_m:\drof{\Sigma V} \to \drof{\Sigma V}. \]
The Harrison cohomology of $V$ with coefficients in $V^*$ is defined as the cohomology of the complex $\caq{V}{V^*}$ and denoted by $\haq{V}{V^*}$.
\end{defi}

We now define the appropriate cyclic cohomology theory; it is closely related to symplectic structures as we will see later on.

\begin{defi}
Let $V$ be a free graded module. Let $m:\clalg{V} \to \clalg{V}$ be a $\ci$-structure. The cyclic Harrison complex of the $\ci$-algebra $V$ is defined on the module consisting of all 0-forms:
\[ \ccaq{V}:=\Sigma\drzf{\Sigma V}. \]
The differential on this complex is the (suspension of the) Lie operator of the vector field $m$;
\[ L_m:\drzf{\Sigma V} \to \drzf{\Sigma V}. \]
The cyclic Harrison cohomology of $V$ is defined as the cohomology of the complex $\ccaq{V}$ and is denoted by $\hcaq{V}$.
\end{defi}

The main result on the cyclic versus non-cyclic Harrison cohomology of $C_\infty$-algebras is that they essentially coincide, as proved in \cite{HL}. We will use the following version of this result  contained in the cited reference, see also \cite{loday}. Note that the noncommutative de Rham differential induces a map $I:\ccaq{V}\rightarrow \caq{V}{V^*}$. If $V$ is a strictly graded commutative algebra then both $\ccaq{V}$ and $\caq{V}{V^*}$ are bigraded; we say an element has homogeneous bidegree $(i,j)$ if it has order $i$ and degree $j$, where a constant 1-form or function is considered to have order zero.

\begin{theorem} \label{cor_caqiso}
Let $V$ be a unital strictly graded commutative algebra. The map
\[ I:\hcaq[i+1,j]{V} \to \haq[ij]{V}{V^*} \]
is;
\begin{enumerate}
\item[(i)]
a monomorphism if $i=1$,
\item[(ii)]
an epimorphism if $i=2$,
\item[(iii)]
an isomorphism if $i \geq 3$.
\end{enumerate}
\end{theorem}
\noproof

\subsubsection{Normalised theories}

We will now provide the analogues of the definitions from the previous subsection for the case of \emph{unital} $\ci$-algebras. Let $V$ be a unital $\ci$-algebra and let $\tau,\boldsymbol{t}$ be a topological basis of $\Sigma V^*$ where $\tau$ is dual to the unit $1 \in V$. We say a vector field $\xi\in\Der(\clalg{V})$ is normalised if it has the form
\[ m = A(\boldsymbol{t})\partial_\tau+\sum_{i\in I}B_i(\boldsymbol{t})\partial_{t_i}. \]
We will denote the Lie subalgebra of $\Der(\clalg{V})$ consisting of all normalised vector fields by $\overline{\Der}(\clalg{V})$.

\begin{defi}
Let $V$ be a unital $\ci$-algebra. The subspace of $\caq{V}{V}$ given by
\[\cnaq{V}{V}:=\Sigma^{-1}\overline{\Der}(\clalg{V})\]
forms a subcomplex of $\caq{V}{V}$ which we refer to as the subcomplex of normalised vector fields.
\end{defi}

Now consider the module of 1-forms $\drof{\Sigma V}$. We say a 1-form $\alpha\in\drof{\Sigma V}$ is normalised if it has the form
\[ \alpha=A(\boldsymbol{t})\cotimes d\tau+\sum_{i\in I}B_i(\boldsymbol{t})\cotimes dt_i. \]
We will denote the subspace of $\drof{\Sigma V}$ consisting of all normalised 1-forms by $\drnof{\Sigma V}$.

\begin{defi}
Let $V$ be a unital $\ci$-algebra. The subspace of $\caq{V}{V^*}$ given by
\[ \cnaq{V}{V^*}:=\Sigma\drnof{\Sigma V}. \]
forms a subcomplex of $\caq{V}{V^*}$ which we refer to as the subcomplex of normalised Lie 1-forms.
\end{defi}

Finally, consider the module of 0-forms $\drzf{\Sigma V}$. We say a 0-form $\alpha\in\drzf{\Sigma V}$ is normalised if it has the form
\[ \alpha = \sum_{i\in I} B_i(\boldsymbol{t})\cotimes t_i. \]
We will denote the subspace of $\drzf{\Sigma V}$ consisting of all normalised 0-forms by $\drnzf{\Sigma V}$.

\begin{defi}
Let $V$ be a unital $\ci$-algebra. The subspace of $\ccaq{V}$ given by
\[ \ccnaq{V} := \Sigma\left[\drnzf{\Sigma V}\right] \]
forms a subcomplex of $\ccaq{V}$, which we refer to as the subcomplex of normalised Lie 0-forms.
\end{defi}

\subsection{Symplectic $\ci$-algebras}

We will now recall (cf. \cite{igusa} and \cite{kontfd}) the definition of a $\ci$-algebra with an invariant inner product. We will call it a \emph{symplectic} $\ci$-algebra in view of its interpretation as a formal symplectic supermanifold together with an odd Hamiltonian vector field. We will also define the appropriate cohomology theory in the symplectic context and show that it essentially reduces to the cyclic theory. Note that the corresponding results also hold for $\ai$-algebras (which we do not consider in the present paper), see \cite{pensch}.

 Recall that $\ci$-structures on a free graded module $V$ can be described in terms of systems of maps
\[ \check{m}_i:V^{\otimes i} \to V, \quad i \geq 1 \]
satisfying certain higher homotopy axioms:

\begin{defi} \label{def_invinn}
Let $V$ be a free graded module of finite rank and let
\[ \innprod{V} \]
be an inner product on $V$.
Let
\[ \check{m}_i:V^{\otimes i} \to V, \quad i \geq 1 \]
be a system of maps determining a $\ci$-structure on $V$. We say that this structure is invariant with respect to the inner product $\langle -,- \rangle$ if the following identity holds for all $x_0,\ldots,x_n \in V$,
\begin{equation} \label{eqn_invinn}
\langle \check{m}_n(x_1,\ldots,x_n),x_0 \rangle = (-1)^{n+|x_0|(|x_1|+\ldots+|x_n|)} \langle \check{m}_n(x_0,\ldots,x_{n-1}),x_n \rangle
\end{equation}
and we say that $V$ is a $\ci$-algebra with an invariant inner product.
\end{defi}

The following result is checked by a straightforward calculation taking place in the completed
tensor algebra of $\Sigma V^*$:

\begin{prop} \label{thm_syminv}
Let $V$ be a free graded module of finite rank. Let
\begin{equation} \label{eqn_syminvdummy}
\check{m}_i:V^{\otimes i} \to V, \quad i\geq 1
\end{equation}
be a $\ci$-structure on $V$ and let
\[ \innprod{V} \]
be an nondegenerate skew-symmetric bilinear form on $V$.

Consider the vector field $m\in \Der(\clalg{V})$
corresponding to the $\ci$-structure and consider the symplectic form $\omega\in \drtf{\Sigma V}$ corresponding to the inner product $\langle -,- \rangle$. The $\ci$-structure \eqref{eqn_syminvdummy} is invariant with respect to the inner product $\langle -,- \rangle$ if and only if $L_m \omega = 0$.
\end{prop}
\noproof

This result motivates the following definition of a \emph{symplectic} $\ci$-algebra:

\begin{defi}
Let $V$ be a free graded module of finite rank. A \emph{symplectic} $\ci$-structure on $V$ is a symplectic form $\omega \in \drtf{\Sigma V}$ together with a \emph{symplectic} vector field
\[ m:\clalg{V} \to \clalg{V} \]
of degree one, such that $m^2=0$.
\end{defi}

\begin{defi}
Let $V$ and $U$ be free graded modules of finite ranks. Let $(m,\omega)$ and $(m',\omega')$ be \emph{symplectic} $\ci$-structures on $V$ and $U$ respectively. A \emph{symplectic} $\ci$-morphism from $V$ to $U$ is a continuous algebra homomorphism
\[ \phi:\clalg{U} \to \clalg{V} \]
of degree zero such that $\phi \circ m' = m \circ \phi$ and $\phi^*(\omega') = \omega$.
\end{defi}

We have the following useful lemmas regarding the cohomology of a symplectic $\ci$-algebra:

\begin{lemma} \label{lem_duaiso}
Let $(V,m,\omega)$ be a \emph{symplectic} $\ci$-algebra. The map $\Phi$ defined by equation \eqref{eqn_ofmvec} is an isomorphism of complexes:
\[\Phi:\caq{V}{V} \to \Sigma^{|\omega|-2}\caq{V}{V^*}.\]
\end{lemma}

\begin{proof}
By Definition \ref{def_sympform} the map $\Phi$ must be bijective. It only remains to prove that this map commutes with the differentials. This will follow from Lemma \ref{lem_schf}. Let $\xi\in\Der(\clalg{V})$ be a vector field, then
\begin{displaymath}
\begin{split}
\Phi(\ad m(\xi)) & = i_{[m,\xi]}(\omega) = [L_m,i_\xi](\omega), \\
& = L_m i_\xi(\omega) = L_m\Phi(\xi). \\
\end{split}
\end{displaymath}
\end{proof}

We need an analogue of Lemma \ref{lem_duaiso} in the context of the cyclic theory:

\begin{lemma} \label{lem_duciso}
Let $(V,m,\omega)$ be a \emph{symplectic} $\ci$-algebra. The Lie subalgebra of
\[\caq{V}{V}:=\left(\Sigma^{-1}\Der(\clalg{V}),\ad m\right)\]
consisting of all \emph{symplectic} vector fields forms a subcomplex denoted by $S\caq{V}{V}$. Furthermore, there is an isomorphism
\[\Upsilon:\Sigma^{|\omega|-2}\ccaq{V} \to S\caq{V}{V}\]
defined by the formula,
\begin{equation} \label{eqn_duciso}
\Upsilon(\alpha):=(-1)^{|\alpha|}\Phi^{-1}d(\alpha);
\end{equation}
where $\Phi$ is the map defined by equation \eqref{eqn_ofmvec}.
\end{lemma}

\begin{proof}
Since the $\ci$-structure $m$ is a symplectic vector field, the Lie subalgebra of symplectic vector fields does indeed form a subcomplex as claimed.

Lemma \ref{lem_schf} part (v) and Lemma \ref{lem_duaiso} tell us that $\Upsilon$ is a map of complexes, i.e. it respects the differentials. It remains to prove that $\Upsilon$ is bijective; this follows from Proposition \ref{prop_symofm} and the (formal) Poincar\'e lemma.
\end{proof}

We will also need one for the normalised cyclic theory:

\begin{lemma} \label{lem_normduciso}
Let $(V,m,\omega)$ be a \emph{unital symplectic} $\ci$-algebra whose symplectic form $\omega$ is \emph{constant} and denote the subcomplex of $\caq{V}{V}$ consisting of all \emph{normalised symplectic} vector fields by $\Scnaq{V}{V}$. The map $\Upsilon$ defined by equation \eqref{eqn_duciso} restricts to a map
\[ \bar{\Upsilon}:\Sigma^{|\omega|-2}\ccnaq{V} \to \Scnaq{V}{V}\]
which is also an isomorphism.
\end{lemma}

\begin{proof}
It is clear that the map $\Phi$ of Lemma \ref{lem_duaiso} restricts to an isomorphism between the subcomplex of normalised vector fields $\cnaq{V}{V}$ and the subcomplex of normalised 1-forms $\cnaq{V}{V^*}$; hence, it remains to show that any closed normalised 1-form is the boundary of a normalised 0-form.

Let $\tau,\boldsymbol{t}$ be a topological basis of $\Sigma V^*$ and let $\Sigma \bar{V}^*$ be the summand generated by the $t_i$. It follows from the results of \cite{HL} that there exists the following commutative diagram:
\[\xymatrix{ \drnzf{\Sigma V} \ar[r]^{d} \ar@<-1ex>@{^{(}->}[d]_{l} & \drnof{\Sigma V} \ar[r]^{d} \ar@<-1ex>@{^{(}->}[d]_{l} & DR^2_{\mathrm{closed}}(\Sigma V) \ar[d]_{\zeta} \\ \prod_{n=2}^\infty [(\Sigma \bar{V}^*)^{\cotimes n}]_{\mathbb{Z}/n\mathbb{Z}} \ar[r]^-{N} \ar@<-1ex>@{->>}[u] & \prod_{n=2}^\infty(\Sigma\bar{V}^*)^{\cotimes(n-1)}\cotimes \Sigma V^* \ar[r]^-{1-z} \ar@<-1ex>@{->>}[u] & \prod_{n=2}^\infty (\Sigma V^*)^{\cotimes n} }\]

The leftmost vertical splitting is given by equation (8.2). The middle vertical splitting is given by equation (7.6). The map $l$ is given by equation (3.5). The map $\zeta$ is defined by Theorem 3.5. The fact that the two vertical splittings commute follows from Corollary 8.2 (i). The map $N$ is the norm map and $z$ is a cyclic permutation. Since the bottom row is clearly exact, the Lemma now follows.
\end{proof}

\section{Finite level structures and Obstruction theory} \label{sec_obstrc}

\subsection{$C_n$-structures}

In this section we define $C_n$-algebras and discuss the problem of lifting a $C_n$-structure to a $C_{n+1}$-structure. The obstruction theory of $n$-algebras is in principle known to the experts, but it is hard to find an explicit reference in the literature, especially in the $C_n$-algebra case. The precursor for this type of obstruction theory is the seminal work of A. Robinson, \cite{robins}, although the main ideas go back to early work on algebraic deformation theory; \cite{gerste}, \cite{nirich}. In the context of unstable homotopy theory, somewhat similar constructions were employed by Stasheff in \cite{stalnm}. We shall only consider the case of \emph{minimal} $C_n$-structures as in this case it is possible to give a convenient interpretation of obstructions in cohomological terms.

Let $V$ be a free graded module. Any vector field $m\in \Der(\clalg{V})$ could be written in the form
\[ m = m_1 + m_2 + \ldots + m_n + \ldots \]
where $m_i$ is a vector field of order $i$. Since we often need to work modulo the endomorphisms of $\clalg{V}$ of order $\geq n$, we will denote the module of such endomorphisms by $(n)$.

\subsubsection{Obstruction theory for $C_n$-algebra structures} \label{sec_obsalg}

In this subsection we will develop the obstruction theory for $C_n$-algebra structures. Let us begin with a definition:

\begin{defi} \label{def_minstr}
Let $V$ be a free graded module. A \emph{minimal} $C_n$-structure ($n \geq 3$) on $V$ is a vector field $m \in \Der(\clalg{V})$ of degree one which has the form
\[ m = m_2 + \ldots + m_{n-1}, \quad m_i \text{ has order $i$} \]
and satisfies the condition $m^2 = 0 \mod (n+1)$.

Let $m$ and $m'$ be two minimal $C_n$-structures on $V$. We say $m$ and $m'$ are equivalent if there is a diffeomorphism $\phi \in \Aut(\clalg{V})$ of the form
\begin{equation} \label{eqn_diffrm}
\phi = \id + \phi_2 + \phi_3 + \ldots + \phi_k +\ldots
\end{equation}
where $\phi_i$ is an endomorphism of order $i$, such that
\[ \phi \circ m \circ \phi^{-1} = m' \mod (n). \]
\end{defi}

\begin{rem}
Recall that we refer to diffeomorphisms of the form \eqref{eqn_diffrm} as \emph{pointed diffeomorphisms}.
\end{rem}

\begin{rem}
Note that this definition is slightly at odds with the definition of $n$-algebras (specifically $A_n$-algebras) given by Stasheff in \cite{staha1} and \cite{staha2}. Every \emph{minimal} $C_n$-algebra is a $C_n$-algebra under Stasheff's definition, however, given a \emph{minimal} $C_n$-structure $m=m_2+\ldots+m_{n-1}$ and an arbitrary vector field $m_n$ of order $n$, $m':=m+m_n$ is a $C_n$-algebra under Stasheff's definition which is obviously not a \emph{minimal} $C_n$-algebra. This distinction will be necessary in order to develop an obstruction theory.
\end{rem}

\begin{rem} \label{rem_ordgde}
Clearly if $m=m_2+\ldots+m_{n-1}$ is a minimal $C_n$-algebra then $m_2$ determines the structure of a strictly graded commutative algebra on the underlying free graded module $V$ which we will call the underlying commutative algebra. Observe that two equivalent minimal $C_n$-structures have the same underlying commutative algebra.

Since the underlying algebra $A:=(V,m_2)$ is a strictly graded commutative algebra, \caq{A}{A} can be given a bigrading: we say a vector field $\xi \in \caq{A}{A}$ has bidegree $(i,j)$ if it is a vector field of order $i$ and has degree $j$ as an element in the profinite graded module $\Sigma^{-1}\Der(\clalg{V})$. The differential on \caq{A}{A} then has bidegree $(1,1)$.
\end{rem}

It will be useful to introduce the following definition:

\begin{defi}
Let $A:=(V,\mu_2)$ be a strictly graded commutative algebra, then the moduli space of minimal $C_n$-structures on $V$ fixing $\mu_2$ is denoted by \mstr{C}{n}{A} and defined as the quotient of the set
\[ \{ m:\clalg{V} \to \clalg{V} : m \text{ is a minimal } C_n\text{-structure and } m_2=\mu_2 \} \]
by the equivalence relation defined in Definition \ref{def_minstr}.
\end{defi}

We will now describe the appropriate terminology which is necessary in discussing extensions of $C_n$-structures to structures of higher order:

\begin{defi}
Let $V$ be a free graded module and let $m = m_2 + \ldots + m_{n-1}$ be a minimal $C_n$-structure on $V$. We say that $m$ is an extendable $C_n$-structure if there exists a vector field $m_n \in \Der(\clalg{V})$ of order $n$ and degree one such that $m + m_n$ is a $C_{n+1}$-structure on $V$ and we call $m_n$ an extension of $m$.

Let $m_n$ and $m'_n$ be two extensions of $m$. We say that $m_n$ and $m'_n$ are equivalent if there exists a diffeomorphism $\phi \in \Aut(\clalg{V})$ of the form
\[ \phi = \id + \phi_{n-1} + \phi_n + \ldots + \phi_{n+k} + \ldots \]
where $\phi_i$ is an endomorphism of order $i$, such that
\begin{equation} \label{eqn_exteqv}
\phi \circ (m + m_n) \circ \phi^{-1} = m + m'_n \mod (n+1).
\end{equation}
The quotient of the set of all extensions of $m$ by this equivalence relation will be denoted by \mext{n}{m}.
\end{defi}

\begin{rem} \label{rem_exteqv}
Let $m_n$ and $m'_n$ be two equivalent extensions of $m$, then there is a diffeomorphism $\phi$ satisfying \eqref{eqn_exteqv} and a vector field $\xi$ of order $n-1$ and degree zero such that $\phi$ has the form
\[ \phi = \id + \xi + \text{ endomorphisms of order}\geq n. \]
From this we conclude that $m'_n = m_n + [\xi,m_2]$.

Conversely if there exists a vector field $\xi$ of order $n-1$ and degree zero such that
\[ m'_n = m_n + [\xi,m_2] \]
then the diffeomorphism $\phi:=\exp(\xi)$ satisfies \eqref{eqn_exteqv}. This means that two extensions of $m$ are equivalent if and only if their difference is a Harrison coboundary.
\end{rem}

\begin{defi} \label{def_obstrc}
Let $V$ be a free graded module and let $m = m_2+\ldots + m_{n-1}$ be a $C_n$-structure on $V$. We define the vector field $\Obs(m)$ of order $n+1$ and degree 2 by
\begin{equation} \label{eqn_obstrc}
\Obs(m) := \frac{1}{2}\sum_{\begin{subarray}{c} i+j = n+2 \\ 3\leq i,j \leq n-1 \end{subarray}} [m_i, m_j].
\end{equation}
\end{defi}

We will now formulate and prove a result analogous to theorems 7.14 and 7.15 of \cite{hamilt}:

\begin{theorem} \label{thm_obstrc}
Let $A:=(V,m_2)$ be a strictly graded commutative algebra. For all $n \geq 3$ \eqref{eqn_obstrc} induces a map
\begin{displaymath}
\begin{array}{ccc}
\mstr{C}{n}{A} & \to & \haq[n+1,3]{A}{A}, \\
m & \mapsto & \Obs(m); \\
\end{array}
\end{displaymath}
which we will denote by $\Obs_n$. The kernel of this map consists of precisely those $C_n$-structures which are extendable:
\[ \{ m \in \mstr{C}{n}{A} : m \textnormal{ is extendable} \} = \ker(\Obs_n). \]
\end{theorem}

\begin{proof}
First of all let us show that given any $C_n$-structure $m=m_2+\ldots+m_{n-1}$ on $V$, $\Obs(m)$ is a cocycle in \caq{A}{A}. Let us define the vector field $\overline{m}$ by the formula
\[ \overline{m}:=m_3+\ldots+m_{n-1} \]
so that we have $m=m_2+\overline{m}$.

By the definition of $\Obs(m)$ (cf. equation \eqref{eqn_obstrc}) we have the equality
\[ m^2 = \Obs(m) \mod (n+2). \]
It follows that
\[ [m_2,\Obs(m)] = \frac{1}{2}[m_2,[m,m]] = [[m_2,m],m] \mod (n+3). \]
Since $A:=(V,m_2)$ is a strictly graded commutative algebra, we know that $[m_2,m_2]=0$. In addition, since $m$ is a $C_n$-structure we know that
\[ [m_2,\overline{m}]+\frac{1}{2}[\overline{m},\overline{m}] = 0 \mod(n+1). \]
From these observations it follows that
\[ [m_2,\Obs(m)] = [[m_2,m],m] = [[m_2,\overline{m}],\overline{m}] = -\frac{1}{2}[[\overline{m},\overline{m}],\overline{m}]=0 \mod(n+3), \]
where the last equality follows from the Jacobi identity, hence $\Obs(m)$ is a cocycle as claimed.


Next we need to show that if
\[ m=m_2+m_3+\ldots+m_{n-1} \quad \text{and} \quad m'=m_2+m'_3+\ldots+m'_{n-1} \]
are two equivalent $C_n$-structures on $V$ then $\Obs(m)$ and $\Obs(m')$ are cohomologous Harrison cocycles.  Since $m$ and $m'$ are equivalent there exists a \emph{pointed} diffeomorphism $\phi \in \Aut(\clalg{V})$ such that
\[ m'=\phi\circ m \circ\phi^{-1} \mod(n). \]
This means there exists a vector field $\xi$ of order $n$ and degree 1 such that
\[ m' = \phi\circ m \circ\phi^{-1} + \xi \mod(n+1). \]

Recall again that by the definition of $\Obs(m)$ we have the equality
\[ m^2 = \Obs(m) \mod (n+2) \]
and likewise we have the same equality for $m'$. Now from the calculation
\begin{displaymath}
\begin{split}
m'^2 & = \phi \circ m^2\circ \phi^{-1} + [m_2,\xi] \mod (n+2), \\
& = \Obs(m) + [m_2,\xi] \mod(n+2); \\
\end{split}
\end{displaymath}
we conclude that $\Obs(m') = \Obs(m) + [m_2,\xi]$.

So far we have proven that the map $\Obs_n$ is well defined. In order to finish the proof we need to show that a $C_n$-structure $m=m_2+\ldots+m_{n-1}$ on $V$ is extendable if and only if $\Obs(m)$ is cohomologous to zero. The reason that this is true is because for any vector field $m_n$ of order $n$ and degree 1 the following identity holds:
\begin{equation} \label{eqn_obstrcdummy}
(m+m_n)^2 = \Obs(m) + [m_2,m_n] \mod(n+2).
\end{equation}
\end{proof}

The next result, also reminiscent of deformation theory, analyses different extensions of a $C_n$-structure in terms of Harrison cohomology:

\begin{theorem} \label{thm_extend}
Let $A:=(V,m_2)$ be a strictly graded commutative algebra and let $m \in \mstr{C}{n}{A}$ be an extendable $C_n$-structure, then \haq[n,2]{A}{A} acts freely and transitively on \mext{n}{m}:
\begin{displaymath}
\begin{array}{ccc}
\haq[n,2]{A}{A} \times \mext{n}{m} & \to & \mext{n}{m}, \\
(\xi_n , m_n) & \mapsto & m_n + \xi_n. \\
\end{array}
\end{displaymath}
\end{theorem}

\begin{proof}
By equation \eqref{eqn_obstrcdummy} we see that $m_n$ is an extension of $m$ if and only if
\begin{equation} \label{eqn_extenddummy}
[m_2,m_n] = -\Obs(m).
\end{equation}
This means that if $m_n$ is an extension of $m$ and $\xi_n$ is a Harrison cocycle then $m_n + \xi_n$ is an extension of $m$. Furthermore if $m_n$ and $m'_n$ are equivalent extensions and $\xi_n$ and $\xi'_n$ are cohomologous cocycles then by Remark \ref{rem_exteqv}, $m_n + \xi_n$ and $m'_n + \xi'_n$ are equivalent extensions of $m$, therefore the above action is well defined.

Condition \eqref{eqn_extenddummy} shows us that if $m_n$ and $m'_n$ are two extensions of $m$ then $m_n - m'_n$ is a Harrison cocycle and hence the above action is transitive. Furthermore if $m_n$ is an extension of $m$ and $\xi_n$ is a Harrison cocycle such that $m_n$ and $m_n + \xi_n$ are equivalent extensions then by Remark \ref{rem_exteqv}, $\xi_n$ is a Harrison coboundary, thus the above action is free.
\end{proof}

\subsubsection{Obstruction theory for $C_n$-algebra morphisms} \label{sec_obsmor}

In this subsection we will develop the obstruction theory for morphisms between two $C_n$-algebras.

\begin{defi} \label{def_minmor}
Let $V$ be a free graded module and let $m$ and $m'$ be two \emph{minimal} $C_n$-structures on $V$. A \emph{minimal} $C_n$-morphism from $m$ to $m'$ is a diffeomorphism $\phi \in \Aut(\clalg{V})$ of degree zero such that
\[ \phi \circ m  = m' \circ \phi \mod (n). \]

Let $\phi$ and $\phi'$ be two such minimal $C_n$-morphisms. We say that $\phi$ and $\phi'$ are homotopic if there exists a vector field $\eta$ of degree $-1$ such that
\[ \phi = \phi'\circ\exp([m,\eta]) \mod (n-1). \]
\end{defi}

\begin{rem}
There is a Lie power series $p(x,y):= x + y + \frac{1}{2}[x,y]+\ldots$ belonging to the pro-free Lie algebra on two generators $x$ and $y$ such that for two vector fields $\gamma$ and $\xi$ of degree zero
\[ \exp(p(x,y)) = \exp(x)\circ\exp(y). \]
It follows that homotopy is an equivalence relation.
\end{rem}

We will now introduce the moduli space of $C_n$-morphisms:

\begin{defi}
Let $V$ be a free graded module and let $m$ and $m'$ be two minimal $C_n$-structures on $V$ with the same underlying algebra, that is to say that $m_2 = m'_2$. The moduli space of minimal $C_n$-morphisms from $m$ to $m'$ is denoted by \mmor{n}{m}{m'} and defined as the quotient of the set
\[ \{ \phi:\clalg{V} \to \clalg{V} : \phi \text{ is a \emph{pointed} minimal $C_n$-morphism from $m$ to $m'$} \} \]
by the homotopy equivalence relation defined in Definition \ref{def_minmor}.
\end{defi}

Let $m = m_2 + \ldots + m_n$ be a $C_{n+1}$-structure on a free graded module $V$. We define the corresponding $C_n$-structure $\bar{m}$ on $V$ as
\[ \bar{m}:= m_2 + \ldots + m_{n-1}. \]
Now we introduce the terminology dealing with extensions of $C_n$-morphisms:

\begin{defi} \label{def_morext}
Let $V$ be a free graded module and let $m$ and $m'$ be two minimal $C_{n+1}$-structures on $V$ with the same underlying algebra. We say that a pointed $C_n$-morphism $\phi$ from $\bar{m}$ to $\bar{m}'$ is extendable if there exists a vector field $\gamma \in \Der(\clalg{V})$ of order $n-1$ and degree zero such that $\exp(\gamma)\circ\phi$ is a $C_{n+1}$-morphism from $m$ to $m'$ and we call $\gamma$ an extension of $\phi$.

Let $\gamma$ and $\gamma'$ be two extensions of $\phi$. We say that $\gamma$ and $\gamma'$ are equivalent if there exists a vector field $\eta$ of order $n-2$ and degree $-1$ such that
\begin{equation} \label{eqn_morextcon}
\exp(\gamma)\circ\phi = \exp(\gamma')\circ\phi\circ\exp([m,\eta]) \mod(n).
\end{equation}
The quotient of the set of all extensions of $\phi$ by this equivalence relation will be denoted by \mmext{n}{\phi}{m}{m'}.
\end{defi}

\begin{rem} \label{rem_morext}
Let $\gamma$ and $\gamma'$ be two extensions of $\phi$ and let $\eta$ be a vector field of order $n-2$ and degree $-1$, then equation \eqref{eqn_morextcon} is satisfied if and only if
\[ \gamma = \gamma' + [m_2,\eta]. \]
It follows that two extensions are equivalent if and only if their difference is a Harrison coboundary in the Harrison complex of the underlying algebra $A:=(V,m_2)$.
\end{rem}

Next we will define the appropriate obstruction to the extension of $C_n$-morphisms:

\begin{defi} \label{def_morobs}
Let $V$ be a free graded module and let $m$ and $m'$ be two minimal $C_{n+1}$-structures on $V$ with the same underlying algebra. Let $\phi$ be a pointed $C_n$-morphism from $\bar{m}$ to $\bar{m}'$. We define the vector field $\obs(\phi)$ of order $n$ and degree $1$ by
\begin{equation} \label{eqn_morobs}
\obs(\phi):= \phi \circ m \circ \phi^{-1} - m' \mod (n+1).
\end{equation}
\end{defi}

We now have all the terminology in place to formulate the analogues of theorems \ref{thm_obstrc} and \ref{thm_extend}:

\begin{theorem} \label{thm_morobs}
Let $A:=(V,m_2)$ be a strictly graded commutative algebra and let $m$ and $m'$ be two minimal $C_{n+1}$-structures ($n \geq 3$) on $V$ whose underlying algebra is $A$, then \eqref{eqn_morobs} induces a map
\begin{displaymath}
\begin{array}{ccc}
\mmor{n}{\bar{m}}{\bar{m}'} & \to & \haq[n,2]{A}{A}, \\
\phi & \mapsto & \obs(\phi); \\
\end{array}
\end{displaymath}
which we will denote by $\obs_n$. The kernel of this map consists of precisely those $C_n$-morphisms which are extendable:
\[ \{ \phi \in \mmor{n}{\bar{m}}{\bar{m}'} : \phi \textnormal{ is extendable} \} = \ker(\obs_n). \]
\end{theorem}

\begin{proof}
Let $\phi$ be a pointed $C_n$-morphism from $\bar{m}$ to $\bar{m}'$, then by Theorem \ref{thm_extend} $\obs(\phi)$ is a cocycle in \caq{A}{A} and $\phi$ is extendable if and only if $\obs(\phi)$ is cohomologous to zero. We need only show that if $\phi'$ is another pointed $C_n$-morphism from $\bar{m}$ to $\bar{m}'$ which is homotopic to $\phi$ then $\obs(\phi)$ and $\obs(\phi')$ are cohomologous cocycles.

Since $\phi$ and $\phi'$ are homotopic there exists a vector field $\gamma$ of order $n-1$ and degree 0 and a vector field $\eta$ of degree $-1$ such that
\[ \phi = \exp(\gamma)\circ\phi'\circ\exp([m,\eta]) \mod (n). \]
Note that for any vector field $\eta$ of degree $-1$,
\[ \exp([m,\eta]) \circ m \circ \exp(-[m,\eta]) = \exp\left(\ad([m,\eta])\right)[m] = m \mod(n+2). \]
We use these facts to demonstrate that the relevant obstructions are cohomologous:
\begin{displaymath}
\begin{split}
\obs(\phi) & = \phi\circ m \circ\phi^{-1}-m' \mod(n+1), \\
& = \exp(\gamma)\circ\phi'\circ\exp([m,\eta]) \circ m \circ \exp(-[m,\eta])\circ\phi'^{-1}\circ\exp(-\gamma) - m' \mod(n+1), \\
& = \phi'\circ m \circ\phi'^{-1} + [\gamma,m_2] - m' \mod(n+1), \\
& = \obs(\phi') + [\gamma,m_2] \mod(n+1). \\
\end{split}
\end{displaymath}
\end{proof}

\begin{theorem} \label{thm_morext}
Let $A:=(V,m_2)$ be a strictly graded commutative algebra and let $m$ and $m'$ be two minimal $C_{n+1}$-structures on $V$ whose underlying algebra is $A$. Let $\phi \in \mmor{n}{\bar{m}}{\bar{m}'}$ be an extendable $C_n$-morphism, then \haq[n-1,1]{A}{A} acts freely and transitively on \mmext{n}{\phi}{m}{m'}:
\begin{displaymath}
\begin{array}{ccc}
\haq[n-1,1]{A}{A} \times \mmext{n}{\phi}{m}{m'} & \to & \mmext{n}{\phi}{m}{m'}, \\
(\xi,\gamma) & \mapsto & \gamma + \xi. \\
\end{array}
\end{displaymath}
\end{theorem}

\begin{proof}
A vector field $\gamma$ of order $n-1$ and degree zero is an extension of $\phi$ if and only if
\begin{equation} \label{eqn_morextdummya}
[m_2,\gamma] = \obs(\phi).
\end{equation}
It follows that if $\gamma$ is an extension of $\phi$ and $\xi$ is a Harrison cocycle then $\gamma + \xi$ is an extension of $\phi$. Moreover if $\gamma$ and $\gamma'$ are equivalent extensions and $\xi$ and $\xi'$ are cohomologous cocycles, then by Remark \ref{rem_morext} $\gamma+\xi$ and $\gamma'+\xi'$ are equivalent extensions; hence the above action is well defined.

If $\gamma$ and $\gamma'$ are two extensions of $\phi$ then by \eqref{eqn_morextdummya} their difference is a cocycle and therefore the above action is transitive. If $\gamma$ is an extension of $\phi$ and $\xi$ is a cocycle such that $\gamma$ and $\gamma+\xi$ are equivalent extensions, then by Remark \ref{rem_morext} $\xi$ is a coboundary and hence the above action is free.
\end{proof}

\subsection{Symplectic $\ci$-structures} \label{sec_symobs}

In this section we develop the obstruction theory in the symplectic context which is parallel to our treatment in section \ref{sec_obstrc}.

In this section we shall only consider \emph{symplectic} $C_n$-algebras with \emph{constant} symplectic forms. This is not really a significant restriction however since the Darboux theorem , cf. \cite{HL} tells us that any symplectic $C_n$-algebra is isomorphic to one of this form. It follows from Proposition \ref{thm_syminv} that a strictly graded commutative \emph{symplectic} $\ci$-algebra $A:=(V,m_2,\omega)$ is the same thing as a commutative Frobenius algebra and we hereafter refer to them as such.

Recall that given a \emph{symplectic} $\ci$-algebra $A:=(V,m,\omega)$, Lemma \ref{lem_duciso} gives us an isomorphism $\Upsilon$ between the (shifted) cyclic Harrison complex $\Sigma^{|\omega|-2}\ccaq{A}$ and the subcomplex of $\caq{A}{A}$ consisting of all symplectic vector fields which we denoted by $S\caq{A}{A}$. Furthermore if $A$ is in fact a strictly graded commutative Frobenius algebra then $\Upsilon$ respects the natural bigrading induced on the cohomology.
\[ \Upsilon: \ccaq[n+1,j+|\omega|-2]{A} \to S\caq[nj]{A}{A}. \]

Recall that given a free graded module $V$, the module consisting of all continuous endomorphisms
\[ f:\clalg{V}\to\clalg{V} \]
of order $\geq n$ is denoted by $(n)$.

\subsubsection{Obstruction theory for symplectic $C_n$-algebra structures}

In this section we will formulate the \emph{symplectic} analogues of the theorems and definitions given in section \ref{sec_obsalg}.

\begin{defi} \label{def_minsym}
Let $V$ be a free graded module of finite rank. A \emph{minimal symplectic} $C_n$-structure ($n \geq 3$) on $V$ consists of a \emph{constant} symplectic form $\omega \in \drtf{\Sigma V}$ and a \emph{symplectic} vector field $m \in \Der(\clalg{V})$ of degree one which has the form
\[ m = m_2 + \ldots + m_{n-1}, \quad m_i \text{ has order $i$} \]
and satisfies the condition $m^2 = 0 \mod (n+1)$.

Let $m$ and $m'$ be two minimal symplectic $C_n$-structures on $V$. We say $m$ and $m'$ are equivalent if there is a \emph{symplectomorphism} of degree zero $\phi \in \Aut(\clalg{V})$ of the form
\begin{equation} \label{eqn_sdifrm}
\phi = \id + \phi_2 + \phi_3 + \ldots + \phi_k + \ldots
\end{equation}
where $\phi_i$ is an endomorphism of order $i$ and such that
\[ \phi \circ m \circ \phi^{-1} = m' \mod (n). \]
\end{defi}

\begin{rem}
Following the convention made in section \ref{sec_notcon} we shall call a symplectomorphism of the form \eqref{eqn_sdifrm} a \emph{pointed symplectomorphism}. Note that under the definition of equivalence, it would be impossible to have two equivalent minimal symplectic $C_n$-structures on $V$ which had different symplectic forms, hence the reason for the omission of the symplectic forms in the definition.
\end{rem}

\begin{rem}
Obviously if $m=m_2+\ldots+m_{n-1}$ is a minimal symplectic $C_n$-algebra with symplectic form $\omega$ then $A:=(V,m_2,\omega)$ is a strictly graded commutative Frobenius algebra which we will call the underlying Frobenius algebra. Observe that two equivalent minimal symplectic $C_n$-structures have the same underlying Frobenius algebra.
\end{rem}

It will be useful to introduce the following definition:

\begin{defi}
Let $A:=(V,\mu_2,\omega)$ be a strictly graded commutative Frobenius algebra, then the moduli space of minimal \emph{symplectic} $C_n$-structures on $V$ fixing $\mu_2$ and $\omega$ is denoted by \mstr[S]{C}{n}{A} and defined as the quotient of the set
\[ \{ m:\clalg{V} \to \clalg{V} : m \text{ is a minimal \emph{symplectic} $C_n$-structure with respect to $\omega$ and $m_2=\mu_2$} \} \]
by the equivalence relation defined in Definition \ref{def_minsym}.
\end{defi}

We will now describe the appropriate terminology necessary in discussing extensions of symplectic $C_n$-structures to structures of higher order:

\begin{defi}
Let $V$ be a free graded module of finite rank and let $\omega \in \drtf{\Sigma V}$ be a constant symplectic form. We say that a minimal symplectic $C_n$-structure $m = m_2 + \ldots + m_{n-1}$ on $V$ is extendable if there exists a \emph{symplectic} vector field $m_n \in \Der(\clalg{V})$ of order $n$ such that $m + m_n$ is a symplectic $C_{n+1}$-structure on $V$ and we call $m_n$ an extension of $m$.

Let $m_n$ and $m'_n$ be two extensions of $m$. We say that $m_n$ and $m'_n$ are equivalent if there exists a \emph{symplectomorphism} $\phi \in \Aut(\clalg{V})$ of degree zero of the form
\[ \phi = \id + \phi_{n-1} + \phi_n + \ldots + \phi_{n+k} + \ldots \]
where $\phi_i$ is an endomorphism of order $i$ and such that
\[ \phi \circ (m + m_n) \circ \phi^{-1} = m + m'_n \mod (n+1). \]
The quotient of the set of all extensions of $m$ by this equivalence relation will be denoted by $\mext[S]{n}{m}$.
\end{defi}

We now identify the appropriate obstruction to extending \emph{symplectic} $C_n$-structures. We will use the same notation as in Definition \ref{def_obstrc}. We will later justify this abuse of notation by Lemma \ref{lem_obscom}:

\begin{defi}
Let $V$ be a free graded module of finite rank and let $\omega \in \drtf{\Sigma V}$ be a constant symplectic form. Let $m = m_2+\ldots + m_{n-1}$ be a minimal symplectic $C_n$-structure on $V$. We define the 0-form $\Obs(m)$ of order $n+2$ and degree $2+|\omega|$ by
\begin{equation} \label{eqn_symobs}
\Obs(m) := \frac{1}{2}\Upsilon^{-1}\left(\sum_{\begin{subarray}{c} i+j = n+2 \\ 3\leq i,j \leq n-1 \end{subarray}} [m_i,m_j] \right).
\end{equation}
\end{defi}

We can now formulate the symplectic analogues of theorems \ref{thm_obstrc} and \ref{thm_extend}. We shall omit the proofs since they are the same as the proofs of section \ref{sec_obstrc} verbatim except that we must use the fact from Lemma \ref{lem_duciso} that $\Upsilon$ is an isomorphism.

\begin{theorem} \label{thm_symobs}
Let $A:=(V,m_2,\omega)$ be a strictly graded commutative Frobenius algebra. For all $n \geq 3$ \eqref{eqn_symobs} induces a map
\begin{displaymath}
\begin{array}{ccc}
\mstr[S]{C}{n}{A} & \to & \hcaq[n+2,1+|\omega|]{A}, \\
m & \mapsto & \Obs(m); \\
\end{array}
\end{displaymath}
which we will denote by $\Obs_n$. The kernel of this map consists of precisely those symplectic $C_n$-structures which are extendable:
\[ \{ m \in \mstr[S]{C}{n}{A} : m \textnormal{ is extendable} \} = \ker(\Obs_n). \]
\end{theorem}
\noproof

\begin{theorem} \label{thm_symext}
Let $A:=(V,m_2,\omega)$ be a strictly graded commutative Frobenius algebra and let $m \in \mstr[S]{C}{n}{A}$ be an extendable symplectic $C_n$-structure, then \hcaq[n+1,|\omega|]{A} acts freely and transitively on \mext[S]{n}{m}:
\begin{displaymath}
\begin{array}{ccc}
\hcaq[n+1,|\omega|]{A} \times \mext[S]{n}{m} & \to & \mext[S]{n}{m}, \\
(\alpha , m_n) & \mapsto & m_n + \Upsilon(\alpha). \\
\end{array}
\end{displaymath}
\end{theorem}
\noproof

\subsubsection{Obstruction theory for symplectic $C_n$-algebra morphisms}

In this section we will develop the obstruction theory for morphisms between two \emph{symplectic} $C_n$-algebras. Again, the formulation is entirely analogous to that of section \ref{sec_obsmor}.

\begin{defi} \label{def_symmor}
Let $V$ be a free graded module of finite rank and let $\omega \in \drtf{\Sigma V}$ be a constant symplectic form. Let $m$ and $m'$ be two minimal \emph{symplectic} $C_n$-structures on $V$. A minimal \emph{symplectic} $C_n$-morphism from $m$ to $m'$ is a \emph{symplectomorphism} $\phi \in \Aut(\clalg{V})$ of degree zero such that
\[ \phi \circ m = m'\circ\phi \mod (n). \]

Let $\phi$ and $\phi'$ be two such minimal symplectic $C_n$-morphisms. We say $\phi$ and $\phi'$ are homotopic if there exists a \emph{symplectic} vector field $\eta$ of degree $-1$ such that
\[ \phi = \phi'\circ\exp([m,\eta]) \mod(n-1). \]
\end{defi}

We will now introduce the moduli space of \emph{symplectic} $C_n$-morphisms:

\begin{defi}
Let $V$ be a free graded module of finite rank and let $\omega \in \drtf{\Sigma V}$ be a constant symplectic form. Let $m$ and $m'$ be two minimal symplectic $C_n$-structures on $V$ with the same underlying Frobenius algebra (i.e. $m_2 = m'_2$). The moduli space of minimal \emph{symplectic} $C_n$-morphisms from $m$ to $m'$ is denoted by \mmor[S]{n}{m}{m'} and defined as the quotient of the set
\[ \{ \phi:\clalg{V} \to \clalg{V} : \phi \text{ is a minimal \emph{pointed symplectic} $C_n$-morphism from $m$ to $m'$} \} \]
by the homotopy equivalence relation defined in Definition \ref{def_symmor}.
\end{defi}

Let $m = m_2 + \ldots + m_n$ be a symplectic $C_{n+1}$-structure on a free graded module $V$. Recall that we define the corresponding symplectic $C_n$-structure $\bar{m}$ on $V$ as
\[ \bar{m}:= m_2 + \ldots + m_{n-1}. \]
Now we will introduce the terminology dealing with extensions of symplectic $C_n$-morphisms:

\begin{defi}
Let $V$ be a free graded module of finite rank and let $\omega \in \drtf{\Sigma V}$ be a constant symplectic form. Let $m$ and $m'$ be two minimal symplectic $C_{n+1}$-structures on $V$ with the same underlying Frobenius algebra. We say that a pointed symplectic $C_n$-morphism $\phi$ from $\bar{m}$ to $\bar{m}'$ is extendable if there exists a \emph{symplectic} vector field $\gamma$ of order $n-1$ and degree zero such that $\exp(\gamma)\circ\phi$ is a symplectic $C_{n+1}$-morphism from $m$ to $m'$ and we call $\gamma$ an extension of $\phi$.

Let $\gamma$ and $\gamma'$ be two extensions of $\phi$. We say that $\gamma$ and $\gamma'$ are equivalent if there exists a \emph{symplectic} vector field $\eta$ of order $n-2$ and degree $-1$ such that
\[ \exp(\gamma)\circ\phi = \exp(\gamma')\circ\phi\circ\exp([m,\eta]) \mod(n). \]
The quotient of the set of all extensions of $\phi$ by this equivalence relation will be denoted by $\mmext[S]{n}{\phi}{m}{m'}$.
\end{defi}

Next we will define the appropriate obstruction to the extension of \emph{symplectic} $C_n$-morphisms. We shall use the same notation as in Definition \ref{def_morobs}. This abuse of notation will later be justified by Lemma \ref{lem_mobscm}.

\begin{defi}
Let $V$ be a free graded module of finite rank and let $\omega \in \drtf{\Sigma V}$ be a constant symplectic form. Let $m$ and $m'$ be two minimal symplectic $C_{n+1}$-structures on $V$ with the same underlying Frobenius algebra. Let $\phi$ be a pointed symplectic $C_n$-morphism from $\bar{m}$ to $\bar{m}'$. We define the 0-form $\obs(\phi)$ of order $n+1$ and degree $1+|\omega|$ by
\begin{equation} \label{eqn_smorob}
\obs(\phi):= \Upsilon^{-1}\left(\phi \circ m \circ \phi^{-1} - m'\right) \mod (n+2).
\end{equation}
\end{defi}

We will now prove the symplectic analogues of theorems \ref{thm_morobs} and \ref{thm_morext}. We shall omit the proofs since they are the same as those of section \ref{sec_obstrc} verbatim except that we must use the fact from Lemma \ref{lem_duciso} that $\Upsilon$ is an isomorphism.

\begin{theorem} \label{thm_smorob}
Let $A:=(V,m_2,\omega)$ be a strictly graded commutative Frobenius algebra and let $m$ and $m'$ be two minimal symplectic $C_{n+1}$-structures ($n \geq 3)$ on $V$ whose underlying Frobenius algebra is $A$, then \eqref{eqn_smorob} induces a map
\begin{displaymath}
\begin{array}{ccc}
\mmor[S]{n}{\bar{m}}{\bar{m}'} & \to & \hcaq[n+1,|\omega|]{A}, \\
\phi & \mapsto & \obs(\phi); \\
\end{array}
\end{displaymath}
which we will denote by $\obs_n$. The kernel of this map consists of precisely those symplectic $C_n$-morphisms which are extendable:
\[ \{ \phi \in \mmor[S]{n}{\bar{m}}{\bar{m}'} : \phi \textnormal{ is extendable} \} = \ker(\obs_n). \]
\end{theorem}
\noproof

\begin{theorem} \label{thm_smorex}
Let $A:=(V,m_2,\omega)$ be a strictly graded commutative Frobenius algebra and let $m$ and $m'$ be two minimal symplectic $C_{n+1}$-structures on $V$ whose underlying Frobenius algebra is $A$. Let $\phi \in \mmor[S]{n}{\bar{m}}{\bar{m}'}$ be an extendable symplectic $C_n$-morphism, then \hcaq[n,|\omega|-1]{A} acts freely and transitively on \mmext[S]{n}{\phi}{m}{m'}:
\begin{displaymath}
\begin{array}{ccc}
\hcaq[n,|\omega|-1]{A} \times \mmext[S]{n}{\phi}{m}{m'} & \to & \mmext[S]{n}{\phi}{m}{m'}, \\
(\alpha,\gamma) & \mapsto & \gamma + \Upsilon(\alpha). \\
\end{array}
\end{displaymath}
\end{theorem}
\noproof

\subsection{Unital Symplectic $\ci$-structures} \label{sec_unisymobs}

In this section we will develop the obstruction theory for \emph{unital symplectic} $\ci$-structures. Here it is our intention to be as brief as possible as the theory is essentially the same as the theories described in the previous two sections. Again, we shall only consider unital symplectic $C_n$-algebras with constant symplectic forms.

Given a \emph{unital} strictly graded commutative Frobenius algebra $A:=(V,m_2,\omega)$, Lemma \ref{lem_normduciso} gives us an isomorphism
\[ \bar{\Upsilon}: \ccnaq[n+1,j+|\omega|-2]{A} \to \Scnaq[nj]{A}{A} \]
between the (shifted) normalised cyclic Harrison complex and the subcomplex of $\caq{A}{A}$ consisting of all normalised symplectic vector fields.

\subsubsection{Obstruction theory for unital symplectic $C_n$-algebra structures}

\begin{defi} \label{def_uniminsym}
Let $V$ be a free graded module of finite rank with basis $1,x_1,\ldots,x_d$. A \emph{minimal unital symplectic} $C_n$-structure ($n \geq 3$) on $V$ consists of a \emph{constant} symplectic form $\omega \in \drtf{\Sigma V}$ and a \emph{symplectic} vector field $m \in \Der(\clalg{V})$ of degree one which has the form
\[ m = m_2 + \ldots + m_{n-1}, \quad m_i \text{ has order $i$} \]
and satisfies the condition $m^2 = 0 \mod (n+1)$.

Furthermore, we impose the condition that for all $i\geq 3$, $m_i$ is a \emph{normalised} vector field and that $m_2$ has the form
\[m_2= \ad\tau-\frac{1}{2}[\tau,\tau]\partial_\tau+\sum_{1\leq i,j,k\leq d}a^k_{ij}[t_i,t_j]\partial_{t_k},\quad a^k_{ij}\in\gf.\]

Let $m$ and $m'$ be two minimal unital symplectic $C_n$-structures on $V$. We say $m$ and $m'$ are equivalent if there is a \emph{unital symplectomorphism} of degree zero $\phi \in \Aut(\clalg{V})$ of the form
\[ \phi = \id + \phi_2 + \phi_3 + \ldots + \phi_k + \ldots\]
where $\phi_i$ is an endomorphism of order $i$ and such that
\[ \phi \circ m \circ \phi^{-1} = m' \mod (n). \]
\end{defi}

It will be useful to introduce the following definition:

\begin{defi}
Let $A:=(V,\mu_2,\omega)$ be a unital strictly graded commutative Frobenius algebra, then the moduli space of minimal \emph{unital symplectic} $C_n$-structures on $V$ fixing $\mu_2$ and $\omega$ is denoted by \mstr[US]{C}{n}{A} and defined as the quotient of the set
\[ \{ m:\clalg{V} \to \clalg{V} : m \text{ is a minimal \emph{unital symplectic} $C_n$-structure and $m_2=\mu_2$} \} \]
by the equivalence relation defined in Definition \ref{def_uniminsym}.
\end{defi}

\begin{defi}
Let $V$ be a free graded module of finite rank and let $\omega \in \drtf{\Sigma V}$ be a constant symplectic form. We say that a minimal unital symplectic $C_n$-structure $m = m_2 + \ldots + m_{n-1}$ on $V$ is extendable if there exists a \emph{normalised symplectic} vector field $m_n \in \overline{\Der}(\clalg{V})$ of order $n$ such that $m + m_n$ is a unital symplectic $C_{n+1}$-structure on $V$ and we call $m_n$ an extension of $m$.

Let $m_n$ and $m'_n$ be two extensions of $m$. We say that $m_n$ and $m'_n$ are equivalent if there exists a \emph{unital symplectomorphism} $\phi \in \Aut(\clalg{V})$ of degree zero of the form
\[ \phi = \id + \phi_{n-1} + \phi_n + \ldots + \phi_{n+k} + \ldots \]
where $\phi_i$ is an endomorphism of order $i$ and such that
\[ \phi \circ (m + m_n) \circ \phi^{-1} = m + m'_n \mod (n+1). \]
The quotient of the set of all extensions of $m$ by this equivalence relation will be denoted by $\mext[US]{n}{m}$.
\end{defi}

\begin{defi}
Let $V$ be a free graded module of finite rank and let $\omega \in \drtf{\Sigma V}$ be a constant symplectic form. Let $m = m_2+\ldots + m_{n-1}$ be a minimal unital symplectic $C_n$-structure on $V$. We define the normalised 0-form $\Obs(m)$ of order $n+2$ and degree $2+|\omega|$ by
\begin{equation} \label{eqn_unisymobs}
\Obs(m) := \frac{1}{2}\bar{\Upsilon}^{-1}\left(\sum_{\begin{subarray}{c} i+j = n+2 \\ 3\leq i,j \leq n-1 \end{subarray}} [m_i,m_j] \right).
\end{equation}
\end{defi}

\begin{theorem} \label{thm_unisymobs}
Let $A:=(V,m_2,\omega)$ be a unital strictly graded commutative Frobenius algebra. For all $n \geq 3$ \eqref{eqn_unisymobs} induces a map
\begin{displaymath}
\begin{array}{ccc}
\mstr[US]{C}{n}{A} & \to & \hcnaq[n+2,1+|\omega|]{A}, \\
m & \mapsto & \Obs(m); \\
\end{array}
\end{displaymath}
which we will denote by $\Obs_n$. The kernel of this map consists of precisely those unital symplectic $C_n$-structures which are extendable:
\[ \{ m \in \mstr[US]{C}{n}{A} : m \textnormal{ is extendable} \} = \ker(\Obs_n). \]
\end{theorem}
\noproof

\begin{theorem} \label{thm_unisymext}
Let $A:=(V,m_2,\omega)$ be a unital strictly graded commutative Frobenius algebra and let $m \in \mstr[US]{C}{n}{A}$ be an extendable unital symplectic $C_n$-structure, then \hcnaq[n+1,|\omega|]{A} acts freely and transitively on \mext[US]{n}{m}:
\begin{displaymath}
\begin{array}{ccc}
\hcnaq[n+1,|\omega|]{A} \times \mext[US]{n}{m} & \to & \mext[US]{n}{m}, \\
(\alpha , m_n) & \mapsto & m_n + \bar{\Upsilon}(\alpha). \\
\end{array}
\end{displaymath}
\end{theorem}
\noproof

\subsubsection{Obstruction theory for unital symplectic $C_n$-algebra morphisms}

\begin{defi} \label{def_unisymmor}
Let $V$ be a free graded module of finite rank with basis $1,x_1,\ldots,x_d$ and let $\omega \in \drtf{\Sigma V}$ be a constant symplectic form. Let $m$ and $m'$ be two minimal \emph{unital symplectic} $C_n$-structures on $V$. A minimal \emph{unital symplectic} $C_n$-morphism from $m$ to $m'$ is a \emph{symplectomorphism} $\phi \in \Aut(\clalg{V})$ of degree zero such that
\[ \phi \circ m = m'\circ\phi \mod (n) \]
and satisfying the condition
\begin{displaymath}
\begin{array}{ccc}
\phi(\tau) & = & \tau + A(\boldsymbol{t}), \\
\phi(t_i) & = & t_i + B_i(\boldsymbol{t}); \\
\end{array}
\end{displaymath}
where $A(\boldsymbol{t}),\{B_i(\boldsymbol{t})\}_{i=1}^d$ are Lie power series in the variables $t_i$.

Let $\phi$ and $\phi'$ be two such minimal unital symplectic $C_n$-morphisms. We say $\phi$ and $\phi'$ are homotopic if there exists a \emph{normalised symplectic} vector field $\eta$ of degree $-1$ such that
\[ \phi = \phi'\circ\exp([m,\eta]) \mod(n-1). \]
\end{defi}

We will now introduce the moduli space of \emph{unital symplectic} $C_n$-morphisms:

\begin{defi}
Let $V$ be a free graded module of finite rank and let $\omega \in \drtf{\Sigma V}$ be a constant symplectic form. Let $m$ and $m'$ be two minimal unital symplectic $C_n$-structures on $V$ with the same underlying Frobenius algebra. The moduli space of minimal \emph{unital symplectic} $C_n$-morphisms from $m$ to $m'$ is denoted by \mmor[US]{n}{m}{m'} and defined as the quotient of the set
\[ \{ \phi:\clalg{V} \to \clalg{V} : \phi \text{ is a minimal \emph{pointed unital symplectic} $C_n$-morphism from $m$ to $m'$} \} \]
by the homotopy equivalence relation defined in Definition \ref{def_unisymmor}.
\end{defi}

Let $m = m_2 + \ldots + m_n$ be a unital symplectic $C_{n+1}$-structure on a free graded module $V$. Recall that we define the corresponding unital symplectic $C_n$-structure $\bar{m}$ on $V$ as
\[ \bar{m}:= m_2 + \ldots + m_{n-1}. \]

\begin{defi}
Let $V$ be a free graded module of finite rank and let $\omega \in \drtf{\Sigma V}$ be a constant symplectic form. Let $m$ and $m'$ be two minimal unital symplectic $C_{n+1}$-structures on $V$ with the same underlying Frobenius algebra. We say that a pointed unital symplectic $C_n$-morphism $\phi$ from $\bar{m}$ to $\bar{m}'$ is extendable if there exists a \emph{normalised symplectic} vector field $\gamma$ of order $n-1$ and degree zero such that $\exp(\gamma)\circ\phi$ is a unital symplectic $C_{n+1}$-morphism from $m$ to $m'$ and we call $\gamma$ an extension of $\phi$.

Let $\gamma$ and $\gamma'$ be two extensions of $\phi$. We say that $\gamma$ and $\gamma'$ are equivalent if there exists a \emph{normalised symplectic} vector field $\eta$ of order $n-2$ and degree $-1$ such that
\[ \exp(\gamma)\circ\phi = \exp(\gamma')\circ\phi\circ\exp([m,\eta]) \mod(n). \]
The quotient of the set of all extensions of $\phi$ by this equivalence relation will be denoted by $\mmext[US]{n}{\phi}{m}{m'}$.
\end{defi}

\begin{defi}
Let $V$ be a free graded module of finite rank and let $\omega \in \drtf{\Sigma V}$ be a constant symplectic form. Let $m$ and $m'$ be two minimal unital symplectic $C_{n+1}$-structures on $V$ with the same underlying Frobenius algebra. Let $\phi$ be a pointed unital symplectic $C_n$-morphism from $\bar{m}$ to $\bar{m}'$. We define the normalised 0-form $\obs(\phi)$ of order $n+1$ and degree $1+|\omega|$ by
\begin{equation} \label{eqn_unismorob}
\obs(\phi):= \bar{\Upsilon}^{-1}\left(\phi \circ m \circ \phi^{-1} - m'\right) \mod (n+2).
\end{equation}
\end{defi}

\begin{theorem} \label{thm_unismorob}
Let $A:=(V,m_2,\omega)$ be a unital strictly graded commutative Frobenius algebra and let $m$ and $m'$ be two minimal unital symplectic $C_{n+1}$-structures ($n \geq 3)$ on $V$ whose underlying Frobenius algebra is $A$, then \eqref{eqn_unismorob} induces a map
\begin{displaymath}
\begin{array}{ccc}
\mmor[US]{n}{\bar{m}}{\bar{m}'} & \to & \hcnaq[n+1,|\omega|]{A}, \\
\phi & \mapsto & \obs(\phi); \\
\end{array}
\end{displaymath}
which we will denote by $\obs_n$. The kernel of this map consists of precisely those unital symplectic $C_n$-morphisms which are extendable:
\[ \{ \phi \in \mmor[US]{n}{\bar{m}}{\bar{m}'} : \phi \textnormal{ is extendable} \} = \ker(\obs_n). \]
\end{theorem}
\noproof

\begin{theorem} \label{thm_unismorex}
Let $A:=(V,m_2,\omega)$ be a unital strictly graded commutative Frobenius algebra and let $m$ and $m'$ be two minimal unital symplectic $C_{n+1}$-structures on $V$ whose underlying Frobenius algebra is $A$. Let $\phi \in \mmor[US]{n}{\bar{m}}{\bar{m}'}$ be an extendable unital symplectic $C_n$-morphism, then $\hcnaq[n,|\omega|-1]{A}$ acts freely and transitively on $\mmext[US]{n}{\phi}{m}{m'}$:
\begin{displaymath}
\begin{array}{ccc}
\hcnaq[n,|\omega|-1]{A} \times \mmext[US]{n}{\phi}{m}{m'} & \to & \mmext[US]{n}{\phi}{m}{m'}, \\
(\alpha,\gamma) & \mapsto & \gamma + \bar{\Upsilon}(\alpha). \\
\end{array}
\end{displaymath}
\end{theorem}
\noproof

\section{The Main Theorem: a Correspondence between $\ci$ and Symplectic $\ci$-structures} \label{correspondence}

In this section we will prove our main result; that a symplectic $\ci$-algebra $(V,m,\omega)$ is uniquely determined by its underlying $\ci$-algebra $(V,m)$ together with the structure of a commutative Frobenius algebra on $(V,m_2)$. The point is that the obstruction theories for symplectic and nonsymplectic $C_\infty$-algebras turn out to be `isomorphic' thanks to Theorem \ref{cor_caqiso}.

\subsection{Main theorem I: the non-unital case}

Once again, in this section we shall only consider symplectic $\ci$-algebras with constant symplectic forms. It will be necessary to introduce the following definitions in order to formulate our main theorem later in this section:

\begin{defi}
\
\begin{enumerate}
\item[(i)]
Let $A:=(V,\mu_2)$ be a strictly graded commutative algebra, then the moduli space of minimal $\ci$-structures on $V$ fixing $\mu_2$ is denoted by $\mstr{C}{\infty}{A}$ and defined as the quotient of the set
\[ \{ m:\clalg{V} \to \clalg{V} : m \text{ is a minimal $\ci$-structure and } m_2=\mu_2 \} \]
by the action under conjugation of the group $G$ consisting of all diffeomorphisms $\phi \in \Aut(\clalg{V})$ of the form
\[ \phi = \id + \phi_2 + \phi_3 + \ldots + \phi_k + \ldots, \]
where $\phi_i$ is an endomorphism of order $i$.
\item[(ii)]
Let $A:=(V,\mu_2,\omega)$ be a strictly graded commutative Frobenius algebra, then the moduli space of minimal \emph{symplectic} $\ci$-structures on $V$ fixing $\mu_2$ and $\omega$ is denoted by \mstr[S]{C}{\infty}{A} and defined as the quotient of the set
\[ \{ m:\clalg{V} \to \clalg{V} : m \text{ is a minimal \emph{symplectic} $\ci$-structure and } m_2=\mu_2 \} \]
by the action under conjugation of the group $G$ consisting of all \emph{symplectomorphisms} $\phi \in \Aut(\clalg{V})$ of the form
\[ \phi = \id + \phi_2 + \phi_3 + \ldots + \phi_k + \ldots, \]
where $\phi_i$ is an endomorphism of order $i$.
\end{enumerate}
\end{defi}

Recall that we call a diffeomorphism (symplectomorphism) of the form
\[ \phi = \id + \phi_2 + \phi_3 + \ldots + \phi_k + \ldots \]
a \emph{pointed} diffeomorphism (symplectomorphism).

\begin{defi}
\
\begin{enumerate}
\item[(i)]
Let $V$ be a free graded module and let $m$ and $m'$ be two minimal $\ci$-structures on $V$ with the same underlying algebra. We say that two $\ci$-morphisms $\phi$ and $\phi'$ from $m$ to $m'$ are homotopic if there exists a vector field $\eta$ of degree $-1$ such that
\[ \phi = \phi'\circ\exp([m,\eta]). \]

We denote the moduli space of $\ci$-morphisms from $m$ to $m'$ by \mmor{\infty}{m}{m'} and define it as the quotient of the set
\[ \{ \phi:\clalg{V} \to \clalg{V} : \phi \text{ is a \emph{pointed} $\ci$-morphism from $m$ to $m'$} \} \]
by the homotopy equivalence relation defined above.
\item[(ii)]
Let $V$ be a free graded module of finite rank and let $\omega \in \drtf{\Sigma V}$ be a constant symplectic form. Let $m$ and $m'$ be two minimal \emph{symplectic} $\ci$-structures on $V$ with the same underlying Frobenius algebra. We say that two \emph{symplectic} $\ci$-morphisms $\phi$ and $\phi'$ from $m$ to $m'$ are homotopic if there exists a \emph{symplectic} vector field $\eta$ of degree $-1$ such that
\[ \phi = \phi'\circ\exp([m,\eta]). \]

We denote the moduli space of \emph{symplectic} $\ci$-morphisms from $m$ to $m'$ by \mmor[S]{\infty}{m}{m'} and define it as the quotient of the set
\[ \{ \phi:\clalg{V} \to \clalg{V} : \phi \text{ is a \emph{pointed symplectic} $\ci$-morphism from $m$ to $m'$} \} \]
by the homotopy equivalence relation defined above.
\end{enumerate}
\end{defi}

Recall that in sections \ref{sec_obstrc} and \ref{sec_symobs} we defined moduli spaces \mstr{C}{n}{A} and \mstr[S]{C}{n}{A}. Given a strictly graded commutative Frobenius algebra $A$ we can define a map
\[ \iota : \mstr[S]{C}{n}{A} \to \mstr{C}{n}{A} \]
which is induced by the canonical inclusion of the Lie subalgebra of symplectic vector fields into the Lie algebra $\Der(\clalg{V})$ of all vector fields. This map is well defined because the group of symplectomorphisms is a subgroup of the group $\Aut(\clalg{V})$ of all diffeomorphisms. This map is also defined for $n=\infty$.

Similarly if $(m,\omega)$ is a symplectic $C_n$-structure on a free graded module $V$ of finite rank then we can define a map
\[ \iota : \mext[S]{n}{m} \to \mext{n}{m} \]
which is again induced by the canonical inclusion of the Lie subalgebra of symplectic vector fields into $\Der(\clalg{V})$.

The same applies to the moduli spaces $\mmor{n}{m}{m'}$ and $\mmor[S]{n}{m}{m'}$ which we also introduced in sections \ref{sec_obstrc} and \ref{sec_symobs}. Given a constant symplectic form $\omega \in \drtf{\Sigma V}$ and two symplectic $C_n$-structures $m$ and $m'$ with the same underlying Frobenius algebra we can define a map
\[ \iota : \mmor[S]{n}{m}{m'} \to \mmor{n}{m}{m'} \]
which is induced by the canonical inclusion of the subgroup of symplectomorphisms into the group $\Aut(\clalg{V})$ of all diffeomorphisms. The fact that this map is well defined modulo the homotopy equivalence relations follows tautologically from the fact that the symplectic vector fields are a Lie subalgebra of the Lie algebra of all vector fields. Likewise this map is also defined for $n=\infty$.

Given a constant symplectic form $\omega \in \drtf{\Sigma V}$, two symplectic $C_{n+1}$-structures $m$ and $m'$ and a symplectic $C_n$-morphism $\phi$ from $\bar{m}$ to $\bar{m}'$ we can also define a map
\[ \iota : \mmext[S]{n}{\phi}{m}{m'} \to \mmext{n}{\phi}{m}{m'} \]
which is induced by the canonical inclusion of the Lie subalgebra of symplectic vector fields into $\Der(\clalg{V})$. Again it is tautological to check that this map is well defined modulo the homotopy equivalence relations.

If $A:=(V,m_2,\omega)$ is a strictly graded commutative Frobenius algebra then recall that there is a natural bigrading on cohomology. The isomorphisms of lemmas \ref{lem_duaiso} and \ref{lem_duciso} respect this bigrading. This gives us the following commutative diagram for all $j \in \mathbb{Z}$ and all $n\geq 1$:
\begin{equation} \label{fig_isocom}
\xymatrix{ \hcaq[n+1,j+|\omega|-2]{A} \ar^{I}[r] \ar^{\Psi}[rd] \ar_{\Upsilon}[d] & \haq[n,j+|\omega|-2]{A}{A^*} \\ S\haq[nj]{A}{A} \ar@{^{(}->}[r] & \haq[nj]{A}{A} \ar^{\Phi}[u] \\ }
\end{equation}
where $\Phi$ is the isomorphism which is defined by equation \eqref{eqn_ofmvec}, $\Upsilon$ is the isomorphism which is defined by equation \eqref{eqn_duciso} and $I$ is the map of Theorem \ref{cor_caqiso} which comes from the periodicity exact sequence. The map $\Psi$ is defined as $\Psi:=\Phi^{-1}\circ I$. Any action of the group $\haq[nj]{A}{A}$ could be pulled back along $\Psi$ to an action of the group $\hcaq[n+1,j+|\omega|-2]{A}$. We will now need the following auxiliary lemmas in order to prove our main result:

\begin{lemma} \label{lem_obscom}
Let $A:=(V,m_2,\omega)$ be a strictly graded commutative Frobenius algebra:
\begin{enumerate}
\item[(i)]
For all $n\geq 3$ the following diagram commutes:
\[ \xymatrix{ \mstr{C}{n}{A} \ar^-{\Obs_n}[r] & \haq[n+1,3]{A}{A} \\ \mstr[S]{C}{n}{A} \ar^{\iota}[u] \ar^-{\Obs_n}[r] & \hcaq[n+2,1+|\omega|]{A} \ar^{\Psi}[u] \\ } \]
\item[(ii)]
Let $m \in \mstr[S]{C}{n}{A}$ be an extendable symplectic $C_n$-structure. The map
\[ \iota : \mext[S]{n}{m} \to \mext{n}{m} \]
is $\hcaq[n+1,|\omega|]{A}$-equivariant.
\end{enumerate}
\end{lemma}

\begin{proof}
This is a tautological consequence of diagram \eqref{fig_isocom}.
\end{proof}

\begin{lemma} \label{lem_mobscm}
Let $A:=(V,m_2,\omega)$ be a strictly graded commutative Frobenius algebra and let $m$ and $m'$ be two minimal symplectic $C_{n+1}$-structures on $V$ whose underlying Frobenius algebra is $A$:
\begin{enumerate}
\item[(i)]
The following diagram commutes:
\[ \xymatrix{ \mmor{n}{\bar{m}}{\bar{m}'} \ar^{\obs_n}[r] & \haq[n,2]{A}{A} \\ \mmor[S]{n}{\bar{m}}{\bar{m}'} \ar^{\iota}[u] \ar^{\obs_n}[r] & \hcaq[n+1,|\omega|]{A} \ar^{\Psi}[u] \\ } \]
\item[(ii)]
Let $\phi \in \mmor[S]{n}{\bar{m}}{\bar{m}'}$ be an extendable symplectic $C_n$-morphism. The map
\[ \iota : \mmext[S]{n}{\phi}{m}{m'} \to \mmext{n}{\phi}{m}{m'} \]
is $\hcaq[n,|\omega|-1]{A}$-equivariant.
\end{enumerate}
\end{lemma}

\begin{proof}
Again this follows tautologically from diagram \eqref{fig_isocom}.
\end{proof}

We are now ready to formulate our main result:

\begin{theorem} \label{thm_main}
Let $A:=(V,m_2,\omega)$ be a strictly graded \emph{unital} commutative Frobenius algebra:
\begin{enumerate}
\item[(i)]
The map
\[ \iota : \mstr[S]{C}{\infty}{A} \to \mstr{C}{\infty}{A} \]
is a bijection.
\item[(ii)]
Let $m$ and $m'$ be two minimal symplectic $\ci$-structures on $V$ whose underlying Frobenius algebra is $A$. The map
\[ \iota : \mmor[S]{\infty}{m}{m'} \to \mmor{\infty}{m}{m'} \]
is a surjection.
\end{enumerate}
\end{theorem}

\begin{proof}
Let us begin by proving that the map $\iota : \mstr[S]{C}{\infty}{A} \to \mstr{C}{\infty}{A}$ is surjective. Let \[ m = m_2 + m_3 + \ldots + m_n + \ldots \]
be a minimal $\ci$-structure on $V$. We will inductively construct a sequence of \emph{symplectic} vector fields $m'_i, \ 3\leq i < \infty$ and a sequence of vector fields $\gamma_i , \ 2\leq i < \infty$, where $m'_i$ has order $i$ and degree one and $\gamma_i$ has order $i$ and degree zero, such that:
\begin{enumerate}
\item[(i)]
\[ m':= m_2 + m'_3 + \ldots + m'_n + \ldots \]
is a minimal \emph{symplectic} $\ci$-structure.
\item[(ii)]
\[ \phi:=\ldots\circ\exp(\gamma_n)\circ\ldots\circ\exp(\gamma_3)\circ\exp(\gamma_2) \]
is a $\ci$-morphism from $m$ to $m'$.
\end{enumerate}

Let us assume that we have constructed a sequence of \emph{symplectic} vector fields $m'_3,\ldots,m'_{n-1}$ of degree one and a sequence of vector fields $\gamma_2,\ldots,\gamma_{n-2}$ of degree zero, where $m'_i$ and $\gamma_i$ have order $i$, satisfying:
\begin{enumerate}
\item[(i)]
\[ m':= m_2 + m'_3 + \ldots + m'_{n-1} \]
is a minimal \emph{symplectic} $C_n$-structure.
\item[(ii)]
\[ \phi:= \exp(\gamma_{n-2})\circ\ldots\circ\exp(\gamma_2) \]
is a minimal $C_n$-morphism from $\bar{m}:= m_2 + \ldots + m_{n-1}$ to $m'$.
\end{enumerate}
Note that the base case $n=3$ is trivial. The $C_n$-structures $\bar{m}$ and $m'$ represent the same class in $\mstr{C}{n}{A}$, therefore by Lemma \ref{lem_obscom} and Theorem \ref{thm_obstrc} we see that
\[ \Psi(\Obs_n(m')) = \Obs_n(\iota(m')) = \Obs_n(\bar{m}) = 0 \]
because the $C_n$-structure $\bar{m}$ is extendable. Theorem \ref{cor_caqiso} now implies that $\Obs_n(m') = 0$ and therefore by Theorem \ref{thm_symobs}, $m'$ is an extendable symplectic $C_n$-structure (recall from diagram \eqref{fig_isocom} that $\Psi$ was defined in terms of the map $I$ of Corollary \ref{cor_caqiso}).

Consider the $\ci$-structure $\phi\circ m\circ\phi^{-1}$; there exists a vector field $\dot{m}_n$ (not necessarily symplectic) of order $n$ and degree one such that
\[ \phi\circ m\circ\phi^{-1} = m' + \dot{m}_n \mod(n+1) \]
and therefore $m'+\dot{m}_n$ is a $C_{n+1}$-structure on $V$. We also know from the above argument that there exists a \emph{symplectic} vector field $m'_n$ of order $n$ and degree one such that $m'+m'_n$ is a \emph{symplectic} $C_{n+1}$-structure. Let us denote the corresponding class of $m'_n$ in $\mext[S]{n}{m'}$ by $[m'_n]$ and the corresponding class of $\dot{m}_n$ in $\mext{n}{m'}$ by $[\dot{m}_n]$. By Theorem \ref{thm_extend} there exists a cohomology class $\xi_n \in \haq[n,2]{A}{A}$ such that
\[ \iota([m'_n]) + \xi_n = [\dot{m}_n]. \]
By Theorem \ref{cor_caqiso} and Lemma \ref{lem_obscom} there exists a cohomology class $\alpha \in \hcaq[n+1,|\omega|]{A}$ such that
\[ \iota([m'_n] + \Upsilon(\alpha)) = [\dot{m}_n]. \]
We see that by modifying our choice of symplectic vector field $m'_n$ appropriately we can assume that there exists a vector field $\gamma_{n-1}$ of order $n-1$ and degree zero such that
\[ \exp(\gamma_{n-1})\circ(m'+\dot{m}_n)\circ\exp(-\gamma_{n-1}) = m' + m'_n \mod(n+1). \]
This completes the inductive step and proves that the map $\iota : \mstr[S]{C}{\infty}{A} \to \mstr{C}{\infty}{A}$ is surjective.

The proof that the map $\iota : \mmor[S]{\infty}{m}{m'} \to \mmor{\infty}{m}{m'}$ is surjective proceeds in a similar fashion to the above proof. We shall outline the main ideas of the proof and leave the reader to provide the details:
\begin{enumerate}
\item[(i)]
Given a $\ci$-morphism $\phi$ from $m$ to $m'$ the idea is to construct by induction, a sequence of \emph{symplectic} vector fields $\gamma'_i, \ 2 \leq i < \infty$ of degree zero and a sequence of vector fields $\eta_i, \ 1 \leq i < \infty$ of degree $-1$, where $\gamma'_i$ and $\eta_i$ have order $i$, such that:
\begin{enumerate}
\item
\[ \phi':=\ldots\circ\exp(\gamma'_n)\circ\ldots\circ\exp(\gamma'_3)\circ\exp(\gamma'_2) \]
is a \emph{symplectic} $\ci$-morphism from $m$ to $m'$.
\item
\[ \phi=\phi'\circ\left[\ldots\circ\exp([m,\eta_n])\circ\ldots\circ\exp([m,\eta_2])\circ\exp([m,\eta_1])\right]. \]
\end{enumerate}
\item[(ii)]
At the $n$th stage use Lemma \ref{lem_mobscm}, Theorem \ref{cor_caqiso} and theorems \ref{thm_morobs} and \ref{thm_smorob} to show that the obstruction to extending the given symplectic $C_n$-morphism to the next level vanishes.
\item[(iii)]
Use Lemma \ref{lem_mobscm}, Theorem \ref{thm_morext} and \emph{both} parts (ii) \emph{and} (iii) of Corollary \ref{cor_caqiso} to modify this extended symplectic $C_{n+1}$-morphism to one which is homotopy equivalent to $\phi$ modulo $(n)$.
\end{enumerate}

Finally we show that the map $\iota : \mstr[S]{C}{\infty}{A} \to \mstr{C}{\infty}{A}$ is injective. Let $m,m'\in\mstr[S]{C}{\infty}{A}$ be two symplectic $\ci$-structures and suppose that there exists a pointed $\ci$-morphism $\phi$ from $m$ to $m'$, then by part (ii) of this theorem there exists a \emph{symplectic} pointed $\ci$-morphism $\phi'$ from $m$ to $m'$ which is homotopy equivalent to $\phi$.
\end{proof}

\begin{rem}
A natural question is whether a result similar to Theorem \ref{thm_main} holds in the $\ai$ and $\li$ contexts. The answer is no. The crucial point on which the proof of Theorem \ref{thm_main} turns is that the cyclic Harrison theory essentially coincides with the (noncyclic) Harrison theory, cf. Theorem \ref{cor_caqiso}. This fails badly for both the associative and Lie cases. For example, let $g$ be a semisimple Lie algebra which could be considered as a symplectic $\li$-algebra together with its Killing form. Then $H^i(g,g)=0$ unless $i=0$; however the cyclic theory $H^\bullet(g,\gf)$ is not zero in dimensions $\geq 3$.

In the associative case a relevant counterexample is provided by the group ring of a finite nonabelian group $G$. Indeed, in this case $\hchoch{\gf[G]}$ is isomorphic to the direct sum of copies of $\hchoch{\gf}\cong \gf[u]$ where the summation is over all conjugacy classes of $G$. On the other hand $\hhoch[i]{\gf[G]}{\gf[G]}=0$ for all $i>0$.
\end{rem}

\begin{rem}
There is a construction, cf. \cite{kontfd} which associates to any minimal symplectic $\ci$-algebra a cycle in an appropriate version of the graph complex. Moreover, two weakly equivalent symplectic $\ci$-algebras give rise to homologous cycles. Therefore, our result shows that the corresponding graph homology class only depends on the underlying $\ci$-algebra together with the structure of a graded commutative Frobenius algebra on its underlying module. In particular, a graph homology class could be associated to any Poincar\'e duality space. It would be interesting to express this construction in classical homotopy theoretic terms, i.e. Massey products.
\end{rem}

\subsection{Main Theorem II: the unital case}

Once again, in this section we shall only consider symplectic $\ci$-algebras with constant symplectic forms.

\begin{defi}
Let $A:=(V,\mu_2,\omega)$ be a unital strictly graded commutative Frobenius algebra, then the moduli space of minimal \emph{unital symplectic} $\ci$-structures on $V$ fixing $\mu_2$ and $\omega$ is denoted by $\mstr[US]{C}{\infty}{A}$ and defined as the quotient of the set
\[ \{ m:\clalg{V} \to \clalg{V} : m \text{ is a minimal \emph{unital symplectic} $\ci$-structure and } m_2=\mu_2 \} \]
by the action under conjugation of the group $G$ consisting of all \emph{symplectomorphisms} $\phi \in \Aut(\clalg{V})$ such that
\begin{displaymath}
\begin{array}{ccc}
\phi(\tau) & = & \tau + A(\boldsymbol{t}), \\
\phi(t_i) & = & t_i + B_i(\boldsymbol{t}); \\
\end{array}
\end{displaymath}
where $A(\boldsymbol{t}),\{B_i(\boldsymbol{t})\}_{i=1}^d$ are Lie power series in the variables $t_i$ consisting of terms of order $\geq 2$. An element of $G$ will be referred to as a \emph{pointed unital symplectomorphism}.
\end{defi}

\begin{defi}
Let $V$ be a free graded module of finite rank and let $\omega \in \drtf{\Sigma V}$ be a constant symplectic form. Let $m$ and $m'$ be two minimal \emph{unital symplectic} $\ci$-structures on $V$ with the same underlying Frobenius algebra. We say that two \emph{unital symplectic} $\ci$-morphisms $\phi$ and $\phi'$ from $m$ to $m'$ are homotopic if there exists a \emph{normalised symplectic} vector field $\eta$ of degree $-1$ such that
\[ \phi = \phi'\circ\exp([m,\eta]). \]

We denote the moduli space of \emph{unital symplectic} $\ci$-morphisms from $m$ to $m'$ by $\mmor[US]{\infty}{m}{m'}$ and define it as the quotient of the set
\[ \{ \phi:\clalg{V} \to \clalg{V} : \phi \text{ is a \emph{pointed unital symplectic} $\ci$-morphism from $m$ to $m'$} \} \]
by the homotopy equivalence relation defined above.
\end{defi}

Given a unital strictly graded commutative Frobenius algebra $A$ we can define a map
\[ \iota : \mstr[US]{C}{\infty}{A} \to \mstr[S]{C}{\infty}{A} \]
which is just defined by forgetting the unital structure.

Given a free graded module $V$ of finite rank, a constant symplectic form $\omega \in \drtf{\Sigma V}$ and two minimal unital symplectic $\ci$-structures $m$ and $m'$ with the same underlying Frobenius algebra, we can define a map
\[ \iota : \mmor[US]{\infty}{m}{m'} \to \mmor[S]{\infty}{m}{m'} \]
which is also defined by forgetting the unital structure.

We are now ready to formulate part II of our main result:

\begin{theorem}
Let $A:=(V,m_2,\omega)$ be a \emph{connected} unital strictly graded commutative Frobenius algebra:
\begin{enumerate}
\item[(i)]
The map
\[ \iota : \mstr[US]{C}{\infty}{A} \to \mstr[S]{C}{\infty}{A} \]
is a bijection.
\item[(ii)]
Let $m$ and $m'$ be two minimal unital symplectic $\ci$-structures on $V$ whose underlying Frobenius algebra is $A$. The map
\[ \iota : \mmor[US]{\infty}{m}{m'} \to \mmor[S]{\infty}{m}{m'} \]
is a surjection.
\end{enumerate}
\end{theorem}

\begin{proof}
The proof proceeds in exactly the same way as the proof of Theorem \ref{thm_main} by making use of the obstruction theory described in sections \ref{sec_symobs} and \ref{sec_unisymobs} and applying the result of \cite{HL} Proposition 8.4 that the map
\[\iota:\hcnaq[ij]{A}\to\hcaq[ij]{A}\]
induced by the inclusion is an isomorphism, except in bidegree $(i,j)=(3,2)$.

A priori, the fact that this map is not surjective in bidegree $(3,2)$ might present a problem in applying an obstruction theory argument; however, since the algebra $A$ is connected, this implies that the degree of the symplectic form $\omega$ is less than or equal to 2. A careful analysis then reveals that none of the cohomology classes relevant to the obstruction theory lie in the critical bidegree $(3,2)$; cf. theorems \ref{thm_unisymobs}, \ref{thm_unisymext}, \ref{thm_unismorob} and \ref{thm_unismorex}.
\end{proof}

\end{document}